\def\version{18 November, 2005}

\documentclass[reqno,11pt]{amsart}

\usepackage{amsmath}
\usepackage{amssymb}

\newfam\Bbbfam
\font\tenBbb=msbm10
\font\sevenBbb=msbm7
\font\fiveBbb=msbm5
\textfont\Bbbfam=\tenBbb
\scriptfont\Bbbfam=\sevenBbb
\scriptscriptfont\Bbbfam=\fiveBbb

\newcommand{\R}     {\mathbb{R}}
\newcommand{\Z}     {\mathbb{Z}}
\newcommand{\N}     {\mathbb{N}}
\renewcommand{\P}   {\mathbb{P}}

\newcommand{\E}     {\mathbb{E}}

\def\1{{\mathchoice {1\mskip-4mu\mathrm l} 
{1\mskip-4mu\mathrm l}
{1\mskip-4.5mu\mathrm l} {1\mskip-5mu\mathrm l}}}
\newcommand{\smallsup}[1] {{\scriptscriptstyle{({#1}})}}
\def\comment#1{}
\newtheoremstyle{thm}{2ex}{2ex}{\itshape\rmfamily}{}
{\bfseries\rmfamily}{}{1.7ex}{}

\newtheoremstyle{rem}{1.3ex}{1.3ex}{\rmfamily}{}
{\itshape\rmfamily}{}{1.5ex}{}

\newenvironment{Proof}[1]
{\vskip0.1cm\noindent{\bf #1. }}{\vspace{0.15cm}}

\newtheorem{theorem}{Theorem}[section]
\newtheorem{lemma}[theorem]{Lemma}
\newtheorem{prop}[theorem] {Proposition}
\newtheorem{cor}[theorem]  {Corollary}

\theoremstyle{definition}
\newtheorem{bem}[theorem] {Remark}
\newtheorem{step}{STEP}

\newcommand{\en}       {\end{equation}}
\newcommand{\eq}       {\begin{equation}}

\newcommand{\eqry}   {\begin{eqnarray}}
\newcommand{\enqry}   {\end{eqnarray}}
\newcommand{\eqarray}   {\begin{eqnarray}}
\newcommand{\enarray}   {\end{eqnarray}}
\newcommand{\eqarraystar} {\begin{eqnarray*}}
\newcommand{\enarraystar} {\end{eqnarray*}}
\newcommand{\bel}{\begin{lemma}}
\newcommand{\el}{\end{lemma}}
\newcommand{\bes}{\begin{step}}
\newcommand{\es}{\end{step}}
\newcommand{\bea}{\begin{array}}
\newcommand{\ea}{\end{array}}
\newcommand{\bpr}{\begin{proof}}
\newcommand{\epr}{\end{proof}}

\renewcommand{\section}{\secdef\sct\sect}
\newcommand{\sct}[2][default]{\refstepcounter{section}
\vspace{0.5cm}
\setcounter{equation}{0}
\centerline{ 
\scshape \arabic{section}.\ #1}
\vspace{0.3cm}}
\newcommand{\sect}[1]{
\vspace{0.5cm}
\centerline{\scshape\Large\bf #1}
\vspace{0.3cm}}

\renewcommand{\subsection}{\secdef \subsct\sbsect}
\newcommand{\subsct}[2][default]{\refstepcounter{subsection}
\nopagebreak
\vspace{0.5\baselineskip}
{\flushleft\bf \arabic{section}.\arabic{subsection}~\bf #1  }
\nopagebreak}
\newcommand{\sbsect}[1]{\vspace{0.1cm}\noindent
{\bf #1}\vspace{0.1cm}}

\renewcommand{\subsubsection}{%
\secdef \subsubsect\sbsbsect}
\newcommand{\subsubsect}[2][default]{%
\refstepcounter{subsubsection}
\nopagebreak
\vspace{0.1\baselineskip}
\nopagebreak
{\flushleft
\sffamily\slshape
\arabic{section}.\arabic{subsection}.\arabic{subsubsection}
\ %
\sffamily #1\/.}\ }
\newcommand{\sbsbsect}[1]{\vspace{0.1cm}\noindent
{\bf #1}\ }

\renewcommand{\d}{{\rm d}}

\newcommand{\eps}{\varepsilon}
\newcommand{\Sym}{\mathfrak{S}}

\newcommand{\ssup}[1] {{\scriptscriptstyle{({#1}})}}
\newcommand{\supp}{{\operatorname {supp}}}
\newcommand{\const}{{\operatorname {cst.}\,}}
\newcommand{\smfrac}[2]{\textstyle{\frac {#1}{#2}}}

\newcommand{\Bcal}   {{\mathcal B }}
\newcommand{\Ccal}   {{\mathcal C }}

\newcommand{\Fcal}   {{\mathcal F }}

\newcommand{\Ical}   {{\mathcal I }}

\newcommand{\Ncal}   {{\mathcal N }}
\newcommand{\Ocal}   {{\mathcal O }}

\newcommand{\Xcal}   {{\mathcal X }}

\begin{document}

\title[Annealed deviations of random walk in random
scenery]{\large Annealed deviations of random walk in random
scenery}
      
\author[Nina Gantert, Wolfgang
        K{\"o}nig and Zhan Shi]{}
\maketitle

\thispagestyle{empty}
\vspace{0.2cm}

\centerline {{\sc Nina Gantert$^1$\/,
Wolfgang  K{\"o}nig$^{2}$\/} and {\sc Zhan Shi$^{3}$\/}}
\vspace{0.8cm}

\centerline {\em $^1$ Fachbereich Mathematik und Informatik der Universit{\"a}t M{\"u}nster,}
\centerline {\em Einsteinstra\ss e 62, D-48149 M{\"u}nster, Germany}
\centerline{\em  
    {\tt gantert@math.uni-muenster.de}}
\vspace{0.3cm}
\centerline{\em $^2$Mathematisches Institut, Universit{\"a}t Leipzig,}
\centerline{\em Augustusplatz 10/11, D-04109 Leipzig, Germany}
\centerline{\tt koenig@math.uni-leipzig.de}
\vspace{0.3cm}
\centerline{\em $^3$ Laboratoire de Probabilit{\'e}s et
    Mod{\`e}les Al{\'e}atoires, Universit{\'e} Paris VI,} 
\centerline{\em 4 place Jussieu, F-75252 Paris Cedex 05, France}
\centerline{\tt  zhan@proba.jussieu.fr}
\vspace{.5cm}

\centerline{\small(\version)}
\vspace{.5cm}

\begin{quote}
  {\small {\bf Abstract:}} Let $(Z_n)_{n\in\N}$ be a
  $d$-dimensional {\it random walk  in  random scenery}, i.e.,
  $Z_n=\sum_{k=0}^{n-1}Y(S_k)$ with $(S_k)_{k\in\N_0}$ a
  random walk in $\Z^d$ and $(Y(z))_{z\in\Z^d}$ an i.i.d.~scenery,
  independent of the walk. The walker's steps have mean
  zero and some finite exponential moments.  
  We identify the speed and the rate of the logarithmic decay of 
  $\P(\frac 1n Z_n>b_n)$ for various choices of sequences $(b_n)_n$ in $[1,\infty)$. Depending 
  on $(b_n)_n$ and the upper tails of the scenery, we identify different regimes for the speed 
  of decay and different variational formulas for the rate functions. In contrast to recent
work \cite{AC02} by A.~Asselah and F.~Castell, we consider sceneries {\it unbounded\/} to infinity.
It turns out that there are interesting connections to large deviation properties of 
self-intersections of the walk, which have been studied recently by X.~Chen \cite{C03}.
\end{quote}

\bigskip

\begin{quote}
  {\small {\bf R\'esum\'e~:}} Soit $(Z_n)_{n\in\N}$ une marche
al\'eatoire en paysage al\'eatoire sur $\Z^d$~; il s'agit du
processus d\'efini par $Z_n=\sum_{k=0}^{n-1}Y(S_k)$, o\`u
$(S_k)_{k\in\N_0}$ est une marche al\'eatoire \`a valeurs dans $\Z^d$, 
et
le paysage al\'eatoire $(Y(z))_{z\in\Z^d}$ est une famille de variables
al\'eatoires i.i.d.\ independante de la marche. On suppose que
$S_1$ est centr\'ee et admet certains moments exponentiels finis.
Nous identifions la vitesse et la fonction de taux de $\P(\frac 1n
Z_n>b_n)$, pour diverses suites $(b_n)_n$ \`a valeurs dans $[1,\infty[$.
Selon le comportement de $(b_n)_n$ et de la queue de distribution du
paysage al\'eatoire, nous d\'ecouvrons diff\'erents r\'egimes ainsi
que diff\'erentes formules variationnelles pour les fonctions de taux.
Contrairement au travail r\'ecent de A.~Asselah and
F.~Castell~\cite{AC02}, nous \'etudions le cas o\`u le paysage 
al\'eatoire
n'est {\it pas born\'e\/}. Finalement, nous observons des liens
int\'eressants avec certaines propri\'et\'es d'auto-intersection de la
marche $(S_k)_{k\in\N_0}$, r\'ecemment \'etudi\'ees par X.~Chen
\cite{C03}.
\end{quote}

\vfill


\eject

\setcounter{section}{0}
\section{Introduction}
 \label{Intro}
\subsection{Model and motivation.}\label{sec-model}

\noindent Let $S=(S_n)_{n\in\N_0}$ be a random walk on $\Z^d$ 
starting at the origin. We denote by $\P$ the underlying
probability measure and by $\E$ the corresponding expectation.
We assume that $\E[S_1]=0$ and $\E[|S_1|^2]<\infty$. 
Defined on the same probability space, let $Y=(Y(z))_{z\in\Z^d}$ be an
i.i.d.~sequence of  random variables, independent of the 
walk. We refer to $Y$ as the {\it random scenery.}
Then the process $(Z_n)_{n\in\N}$ defined by
$$
Z_n=\sum_{k=0}^{n-1} Y(S_k),\qquad n\in\N,
$$
where $\N=\{1, 2, \ldots\},$ is called a {\it random walk in random scenery}, sometimes
also referred to as the {\it Kesten-Spitzer random walk in random scenery}, see \cite{KS79}. An
interpretation is as follows. If a random walker has to pay
$Y(z)$ units at any time he/she visits the site $z$, then $Z_n$
is the total amount he/she pays by time $n-1$. 

The random walk in random scenery has been introduced and
analyzed for dimension $d\not=2$ by H.~Kesten and F.~Spitzer \cite{KS79} and
by E.~Bolthausen \cite{B89} for $d= 2$. The case $d=1$ was treated independently by A. ~N. ~Borodin \cite{Bo79a}, \cite{Bo79b}.  Under the assumption that
$Y(0)$ has expectation zero and variance $\sigma^2\in(0,\infty)$, their results imply that
\begin{equation}
    \label{weakconv}
   \frac 1n Z_n\approx a_n^{\ssup{0}}    \begin{cases}n^{-\frac 14}&\mbox{if }d=1,\\
    (\frac n{\log n})^{-\frac 12}&\mbox{if }d=2,\\
    n^{-\frac 12}&\mbox{if }d\ge 3.
    \end{cases}
\end{equation}
More precisely, 
$\frac 1{na_n^{\ssup{0}}}Z_n$ converges in distribution towards some non-degenerate
random variable. The limit is Gaussian in $d\ge 2$ and a
convex combination of Gaussians (but not Gaussian) in $d=1$. This 
can be roughly explained as follows.
In terms of the so-called {\it local times\/} of the walk,
\begin{equation}
    \label{loctim}
    \ell_n(z)=\sum_{k=0}^{n-1} \1_{\{S_k=z\}},\qquad n\in\N, \;\;
    z\in\Z^d,
\end{equation}

\noindent the random walk in random scenery may be identified as
\begin{equation}
    \label{Znrepr}
    Z_n=\sum_{z\in\Z^d}Y(z)\ell_n(z).
\end{equation}

\noindent The number of effective summands in \eqref{Znrepr} is
equal to the {\it range\/} of the walk, i.e., the number of sites 
visited by time $n-1$. Hence, conditional on the
random walk, $Z_n$ is, for dimension $d\ge 3$, a sum of
$\Ocal(n)$ independent copies of  finite multiples of $Y(0)$,
and hence it is plausible that $Z_n/n^{1/2}$  converges to a
normal variable.  The same assertion with logarithmic 
corrections is also plausible in $d=2$. However,  in $d=1$, $Z_n$ is 
roughly a sum of 
$\Ocal(n^{1/2})$ copies of independent variables with
variances of order $\Ocal( n)$, and this suggests the
normalization in \eqref{weakconv}  as well as a non-normal limit.

In this paper, we  analyse deviations $\{\frac 1n Z_n>b_n\}$ for 
various choices of sequences $(b_n)_{n\in\N}$ in $[1,\infty)$. We 
determine the speed and the rate of the logarithmic
asymptotics of the probability of this event as $n\to\infty$, and we
explain the typical behaviour of the random walk and the
random scenery on this event.

This problem has been addressed in recent work \cite{CP01}, \cite{AC02} and 
\cite{Ca04} by F.~Castell in partial collaboration with F.~Pradeilles and 
A.~Asselah for Brownian motion instead of random walk.
 While \cite{CP01} and \cite{Ca04} treat the case of a continuous Gaussian scenery for $b_n=n^{1/2}$
and $\const\leq b_n\ll n^{1/2}$, respectively, the case of an arbitrary 
bounded scenery (constant on the unit cubes) and $b_n=\const$
is considered in \cite{AC02}.
See also \cite{AC02} for further references on this topic and \cite{AC05a} and \cite{GHK04}
for recent results on the random walk case.

The main novelty of the present paper is the study of arbitrary sceneries {\it unbounded\/} 
to $+\infty$ and general scale functions $b_n \geq \const$ in the discrete setting. 
On the technical side, in particular the proof of the upper bound is rather demanding 
and requires new techniques. We solve this part of the problem by a careful analysis 
of high integer moments, a technique which has been recently established in the study 
of intersection properties of random motions.

A very rough, heuristic explanation of the interplay between the deviations of the 
random walk in random scenery and the tails of the scenery at infinity and the dimension $d$ is as follows. 
In order to realize the event $\{\frac 1n Z_n>b_n\}$, it is clear
that the scenery has to assume larger values on the range of the walk than usual.
In order to keep the probabilistic cost for this low, the random walker
has to keep its range small, i.e., it has to concentrate on less sites by time $n$ than usual. 
The {\it optimal\/} joint strategy of the scenery and the walk is determined by
a balance between the respective costs. The optimal strategies in the cases considered
in the present paper are homogeneous. More precisely, the scenery and the walk
each approximate optimal (rescaled) profiles in a large, $n$-dependent box.
These optimal profiles are determined by a (deterministic) variational problem.

The topic of the present paper has deep connections to large deviation properties of 
self-intersections of the walk. This is immediate in the important special case of a 
standard Gaussian scenery $Y$. Indeed, the conditional distribution of $Z_n$ given
 the random walk $S$ is a centered Gaussian with variance equal to 
\begin{equation}\label{Lambdadef}
\Lambda_n=\sum_{z\in \Z^d}\ell_n(z)^2=\|\ell_n\|_2^2,
\end{equation}
which is often called the {\it self-intersection local time}. Hence, large deviations for the random walk in Gaussian scenery would be a consequence of an appropriate large deviation statement for self-intersection local times. However, the latter problem is notoriously difficult and is, up to the best of our knowledge, open in the precision we would need in the present paper. (However, compare to interesting and deep work on self-intersections and mutual intersections by X.~Chen \cite{C03}.) Recent results for self-intersection local times for random walks in dimension $d \geq 5$  and applications to random walk in random scenery are given in \cite{AC05b}.

The remainder of Section~\ref{Intro} is organized as follows. Our main results are in Section~\ref{results}, a heuristic derivation may be found in Section~\ref{sec-heur}, a partial result for Gaussian sceneries for dimension $d=2$ is in Section~\ref{Sec-Gauss}. The structure of the remainder of the paper is as follows. In Section~\ref{sec-varform} we analyse the variational formulas, in Section~\ref{sec-prep} we present the tools for our proofs of the main results, in Sections~\ref{s:ub} and \ref{s:lb} we give the proofs of the upper and the lower bounds, respectively, and finally in the appendix, Section~\ref{sec-proofLDP}, we provide the proof of a large deviation principle that is needed in the paper.

\subsection{Results}\label{results} 

\noindent Our precise assumptions on the random walk, $S$, are the following. The walker starts at $S_0=0$, and the steps have mean zero and some finite exponential moments, more precisely, 
\begin{equation}
\label{emom}
\E[e^{t |S_1|}]<\infty\quad\mbox{for some }t>0.
\end{equation}
  By $\Gamma\in\R^{d\times d}$ we denote the covariance matrix of the walk's step distribution. Hence,
 $S$ lies in the domain of attraction of the  Brownian motion with covariance matrix $\Gamma$.  We assume that $\Gamma$ is a regular matrix. Furthermore, we assume that $S$ is strongly aperiodic, i.e., for any $z\in\Z^d$, the smallest subgroup of $\Z^d$ that contains $\{z+x\colon \P(S_1=x)>0\}$ is $\Z^d$ itself. Finally, to avoid technical difficulties, we also assume that the transition function of the walk is symmetric, i.e., $p(0,z)= p(0,-z)$ for $z\in\Z^d$, where $p(z,\widetilde z)$ denotes the walker's one-step probability from $z\in\Z^d$ to $\widetilde z\in\Z^d$.

Our assumptions on the scenery are the following. Let $Y=(Y(z))_{z\in\Z^d}$ be a family of i.i.d.~random variables, not necessarily having finite expectation, such that 
\begin{equation}\label{assumption}
\E[e^{t Y(0)}]<\infty\quad\mbox{for every }t>0.
\end{equation}
In particular, the {\it cumulant generating function \/}
of $Y(0)$, is finite:
\begin{equation}\label{Hdef}
H(t)=\log\E[e^{tY(0)}]<\infty,\qquad t>0.
\end{equation}
In some of our results, we additionally suppose the following.

\medskip
\noindent{\bf Assumption (Y).} {\it There are constants $D>0$ and $q>1$ 
such
that}
$$
\log\P(Y(0)>r)\sim -D r^q,\qquad r\to\infty.
$$
According to Kasahara's exponential Tauberian theorem (see \cite[Th.~4.12.7]{BGT87}), Assumption (Y) is equivalent to 
\begin{equation}\label{cumgenfct}
H(t)\sim \widetilde D t^p,\quad \mbox{as }t\to\infty,
\qquad\mbox{where}\qquad
\widetilde D=(q-1)(D q^q)^{1/(1-q)}\qquad\mbox{and}\qquad\frac 1q+\frac 1p=1.
\end{equation}

In our first main result, we consider the case of sequences $(b_n)_n$
tending to infinity slower than $n^{\frac 1q}$. By $\nabla$ we denote the usual gradient acting on sufficiently regular functions
$\R^d\to\R$. By $H^1(\R^d)$ we denote the usual Sobolev space, and we write $\|\nabla \psi\|_2^2=\int_{\R^d}|\nabla\psi(x)|^2\,\d x$. 
We use the notation $b_n\gg c_n$ if $\lim_{n\to\infty}b_n/c_n=\infty$.

\begin{theorem}[Very large deviations]\label{inter} Suppose that 
Assumption (Y) holds with some $q > \frac d2$. Pick a sequence $(b_n)_{n\in\N}$ satisfying $1\ll b_n\ll
n^{\frac 1q}$. Then
\begin{equation}\label{interasy}
\lim_{n\to\infty}n^{-\frac d{d+2}} 
b_n^{-\frac {2q}{d+2}}\log \P(\smfrac 1n Z_n>b_n)=- K_{D,q},
\end{equation}
where 
\begin{equation}\label{kdq}
K_{D,q}\equiv  \inf\Bigl\{
    \frac 12\|\Gamma^{\frac 12}\nabla \psi\|_2^2+D\|\psi^2\|_{p}^{-q}\colon 
    \psi\in H^1(\R^d),\|\psi\|_2=1\Bigr\},
\end{equation}
(we recall that $\frac 1p+\frac 1q=1$), and $K_{D,q}$ is positive.
\end{theorem}

\begin{bem}\label{rem-inter} For $q \in(1,\frac d2)$, (\ref{interasy}) also holds true, but $K_{D,q} = 0$. Indeed, this follows from Proposition~\ref{constpos} below together with our proof of Theorem~\ref{inter}. One can also see this directly by giving an explicit  lower bound for $\log\P(\smfrac 1n Z_n>b_n)$ which runs on a strictly smaller scale than $n^{\frac d{d+2}} 
b_n^{\frac {2q}{d+2}}$. 
It remains an open problem in this paper to determine the {\it precise\/} logarithmic rate of 
$\P(\smfrac 1n Z_n>b_n)$ in the case $q\in(1,\frac d2)$. 
 The case $q=\frac d2$ seems even more delicate and is also left open in the present paper.
The case $q\in(0,1)$ has been studied in \cite{GHK04}.
\hfill$\Diamond$
\end{bem}

Note that the variational problem in \eqref{kdq} is of independent interest;
it also appeared in \cite[Theorem 1.1]{BAL91} in the context of heat
kernel asymptotics.
In Proposition~\ref{constpos} below it turns out that $K_{D,q}$ is positive if and only if $q\geq \frac d2$. 

Our next result essentially extends \cite[Th.~2.2]{AC02} from the case of bounded sceneries to the case in \eqref{assumption}. 

\begin{theorem}[Large deviations]\label{lin} Suppose that \eqref{assumption} holds.
Assume that $\E[Y(0)]=0$, and set $\overline p\equiv\limsup_{t\to\infty}\frac{\log H(t)}{\log t}$. Assume that $\overline p<\infty$ in $d\leq 2$ respectively $\overline p<\frac d{d-2}$ in $d\geq 3$. Then, for any $u>0$ satisfying $u\in \supp(Y(0))^\circ$, 
\begin{equation}\label{linasy}
\lim_{n\to\infty}n^{-\frac d{d+2}}\log\P(\smfrac 1n Z_n>u)= -K_{H}(u),
\end{equation}
where 
\begin{equation} \label{K3def}
K_H(u)\equiv \inf\Bigl\{
    \frac 12\|\Gamma^{\frac 12}\nabla \psi\|_2^2+\Phi_H(\psi^2,u)\colon \psi\in 
    H^1(\R^d),\|\psi\|_2=1\Bigr\},
\end{equation}
and
\begin{equation}\label{Phidef}
\Phi_H(\psi^2,u)=\sup_{\gamma\in(0,\infty)}\Bigl[\gamma 
u-\int_{\R^d}H(\gamma \psi^2(y))\,\d y\Bigr]. 
\end{equation}
The constant $K_H(u)$ is positive.
\end{theorem}

Switching to the scenery $-Y$, one may, under appropriate conditions, use Theorem~\ref{lin} to obtain the \lq other half\rq\ of a full large deviation principle for $(\frac 1n Z_n)_n$. This was carried out in \cite{AC02} for bounded sceneries. For Brownian motion in a Gaussian scenery, a result analogous to Theorems~\ref{inter} and \ref{lin} is \cite[Th.~2]{Ca04}. 

Note that the constant $K_H(u)$ depends on the entire scenery distribution, while $K_{D,q}$ in \eqref{kdq} only depends on its upper tails.

\begin{bem} A statement analogous to Remark~\ref{rem-inter} also applies here: for dimensions $d \geq 3$, when $ \liminf_{t\to\infty}\frac{\log H(t)}{\log t} > \frac d{d-2}$, (\ref{linasy}) also holds true, but $K_H(u) = 0$ for any $u>0$. It was shown recently in \cite{AC05a} that under assumption (Y) with $q\in(1,  \frac{d}{2})$ and an additional symmetry assumption, $\log\P(\smfrac 1n Z_n>b_n)$ is of the order $n^{\frac q{q+2}}$. The case $q\in(0,1)$ has been studied in \cite{GHK04}.
\hfill$\Diamond$
\end{bem}

\begin{bem}[Large deviations and non-convexity]\label{rem-LDP} It is easy to see that, in the special case where $H(t)=\widetilde Dt^p$ (see \eqref{cumgenfct}), $K_H(u)=u^{\frac{2q}{d+2}} K_{D,q}$, for any $u>0$. (For asymptotic scaling relations see Lemma~\ref{constasy}.) In particular, $\frac 1n Z_n$ satisfies a large deviation principle on $(0,\infty)$ with speed $ n^{\frac d{d+2}}$ and rate function $u\mapsto u^{\frac{2q}{d+2}} K_{D,q}$. This function is strictly convex for $q>\frac d2 +1$ and strictly concave for $q<\frac d2 +1$. In the important special case of a centered Gaussian scenery, Theorem~\ref{lin} contains non-trivial information only in the case $d\in\{1,2,3\}$, in which the rate function is strictly convex, linear and strictly concave, respectively; see also \cite{CP01} and \cite{Ca04}.

The non-convexity around zero for bounded sceneries in $d\in\{3,4\}$ was found in \cite{AC02} by proving that $K_H(u)\geq C u^{\frac 4{d+2}}$ as $u\to0$ for some positive constant $C$. 
\hfill$\Diamond$
\end{bem}

The upper bounds in Theorems~\ref{inter} and \ref{lin} are proved in Section~\ref{s:ub}, and the lower bounds in Section~\ref{s:lb}. We consider only sequences $b_n\ge 1$ there. The case $a_n^{\ssup{0}}\ll b_n\ll 1$ seems subtle and is left open in the present paper; however see Section~\ref{Sec-Gauss} for a partial result.

Our next proposition gives almost sharp criteria for the positivity of the constants $K_{D,q}$ and $K_H(u)$ appearing in Theorems~\ref{inter} and \ref{lin}. 

\begin{prop} [Positivity of the constants]\label{constpos} Fix $d\in\N$ and $p,q>1$ satisfying $\frac 1p+\frac 1q=1$.
\begin{enumerate}
\item[(i)] For any $D>0$,
\begin{equation}\label{Kdef}
K_{D,q}= (d+2)\Bigl(\frac D2\Bigr)^{\frac 2{d+2}}\Bigl(\frac{
\chi_{d,p}}d\Bigr)^{\frac d{d+2}},
\end{equation}
where
\begin{equation}\label{chiident}
\chi_{d,p}\inf\Bigl\{\frac 1{2}\|\Gamma^{\frac 12}\nabla \psi\|_2^2\colon \psi\in 
H^1(\R^d)
\colon\|\psi\|_2=1=\|\psi\|_{2p}\Bigr\}.
\end{equation}
The constant $\chi_{d,p}$ is positive if and only if $d\leq \frac {2p}{p-1}=2q$. 
Hence, $K_{D,q}$ is positive if and only if $d\leq \frac {2p}{p-1}=2q$.

\item[(ii)] The constant $K_H(u)$ is positive for any $u>\E[Y(0)]=0$ if 
$$
\limsup_{t\to\infty}\frac{\log H(t)}{\log 
t}<\begin{cases}\infty&\mbox{if }d\leq 2,\\
\frac d{d-2}&\mbox{if }d\geq3.
\end{cases}
$$
For $d\geq 3$, if $\liminf_{t\to\infty}\frac{\log H(t)}{\log t}>\frac d{d-2}$, then $K_H(u)=0$ for any $u>0$.
\end{enumerate}
\end{prop}

The proof of Proposition \ref{constpos} is in Section~\ref{sec-varform}. There we also clarify the relation between $\chi_{d,p}$ and the so-called {\it Gagliardo-Nirenberg\/} constant.

Now we formulate asymptotic relations between the rates obtained in Theorems~\ref{inter} and \ref{lin}.

\begin{lemma}[Asymptotic scaling relations]\label{constasy} Fix $D>0$ and $q>1$, and recall \eqref{cumgenfct}. 
\begin{enumerate}
\item [(i)] Assume that $H(t)\sim\widetilde D t^p$ as $t\to\infty$, then
\begin{equation}\label{Kasy2}
K_H(u)\sim u^{\frac {2q}{d+2}} K_{D,q}\qquad\mbox{as }u\to\infty.
\end{equation}

\item [(ii)] Assume that $\E[Y(0)]=0$ and $\E[Y(0)^2]=1$, then
\begin{equation}\label{Kasy3}
K_H(u)\leq u^{\frac {4}{d+2}} \bigl[K_{\frac 12,2}+o(1)\bigr]\qquad\mbox{as }u\downarrow 0.
\end{equation}

\end{enumerate}
\end{lemma}

The proof of Lemma~\ref{constasy} is in Section~\ref{Sec-ScalLim}. 

\begin{bem} We conjecture that the lower bound in \eqref{Kasy3} also holds under an appropriate upper bound on $H$. It is clear (see Remark~\ref{rem-LDP} and note the monotonicity of $K_H(u)$ in $H$) that $u^{-4/(d+2)}K_H(u)\geq K_{D,2}$ for every $u>0$ if $H(t)\leq \widetilde D t^2$ for every $t\geq 0$. The positivity of $\liminf_{u\downarrow 0} u^{-4/(d+2)}K_H(u)$ (for cumulant generating functions $H$ of {\it bounded\/} variables) is contained in \cite{AC02} as part of the proof for non-convexity of the rate function $K_H$ in $d\in\{3,4\}$. Since $K_{D,2}=0$ in $d> 4$, it is clear that this proof must fail in $d> 4$.
\hfill$\Diamond$
\end{bem}

Lemma~\ref{constasy}(i) is consistent with Theorems~\ref{inter} and \ref{lin}. 

\subsection{Heuristic derivation of Theorems~\ref{inter} and \ref{lin}}\label{sec-heur}

\noindent The asymptotics in \eqref{interasy} and \eqref{linasy}
are based on large deviation principles for scaled versions of the 
walker's local times $\ell_n$ and the scenery $Y$.
A short summary of the joint optimal strategy of the walker and the scenery is the following.
Let us first explain the exponential decay rate of the probabilities under consideration.
Assume that $1\ll b_n\ll n^{\frac 1q}$.
In order to contribute optimally to the event $\{\frac 1n Z_n>b_n\}$, the walker spreads out over a
region whose diameter is of order $\alpha_n$ (for a particular choice of 
$\alpha_n$, depending on the sequence $(b_n)_n$). The cost for this behavior is 
$e^{\Ocal(n\alpha_n^{-2})}$. The scenery assumes extremely
large values within that region, more precisely: values of the order $b_n$.
The cost for doing that is $\exp\{\Ocal(b_n^q\alpha_n^d)\}$, under Assumption (Y).
The choice of $\alpha_n$ is now determined by putting
\begin{equation}\label{alphachoice}
\frac n{\alpha_n^2}=\alpha_n^d b_n^q.
\end{equation}
A calculation shows that for this choice of $\alpha_n$ both sides of
\eqref{alphachoice} are equal to the logarithmic decay order of the probability
$\P(\frac 1n Z_n>b_n)$ in Theorem~\ref{inter}.

Next we give a more precise argument for the very large deviations (Theorem~\ref{inter}) which also explains the constants on the right hand side 
of \eqref{interasy}. Introduce the scaled and normalized
version of the walker's local times, 
\begin{equation}\label{Lndef}
L_n(x)=\frac{\alpha_n^d}n \ell_n\bigl(\lfloor 
x\alpha_n\rfloor\bigr),\qquad x\in \R^d.
\end{equation}
Then $L_n$ is a random element of the set
\begin{equation}
\Fcal=\Bigl\{\psi^2\in L^1(\R^d)\colon \|\psi\|_2=1\Bigr\}
\end{equation}
of all Lebesgue probability densities on $\R^d$. 
Furthermore, introduce the scaled version of the field,
\begin{equation}\label{Yscaled}
\overline Y_n(x)=\frac 1{b_n}Y\bigl(\lfloor 
x\alpha_n\rfloor\bigr),\qquad x\in \R^d.
\end{equation}

Then we have, writing $\langle\cdot,\cdot\rangle$ for the inner product on $L^2(\R^d)$,
\begin{equation}\label{Znscale}
\smfrac 1n Z_n=\smfrac 1n\sum_{z\in\Z^d}\frac n{\alpha_n^d}L_n\bigl({\textstyle{\frac z{\alpha_n}}}\bigr)
{b_n}\overline Y_n\bigl({\textstyle{\frac z{\alpha_n}}}\bigr)=b_n\langle L_n,\overline Y_n\rangle.
\end{equation}
Hence, the logarithmic asymptotics of the probability $\P(\smfrac 1n Z_n>b_n)
=\P(\langle L_n,\overline Y_n\rangle>1)$ will be 
determined by a combination of large deviation principles for $L_n $ and $\overline Y_n$.

In the spirit of the celebrated large deviation
theorem of Donsker and Varadhan, the distributions of $L_n$
satisfy a weak large deviation principle in the weak $L^1$-topology on 
$\Fcal$ with speed $n\alpha_n^{-2}$ and rate function $\Ical \colon 
\Fcal\to[0,\infty]$ given by
\begin{equation}\label{Idef}
\Ical(\psi^2)=\begin{cases} \frac 1{2} \bigl\Vert\Gamma^{\frac 
12}\nabla\psi \bigr\Vert_2^2 
&\mbox{if } \psi\in H^1(\R^d),\\ 
\infty&\mbox{otherwise.} 
\end{cases}
\end{equation}
Roughly speaking, this principle says that, for $\psi^2\in \Fcal$, 
\begin{equation}\label{LDP1heur}
\P(L_n\approx \psi^2)\approx \exp\Big\{-\frac n{\alpha_n^2}\Ical(\psi^2)\Big\},\qquad n\to\infty.
\end{equation}
Using Assumption (Y), we see that the distributions of $\overline Y_n$
should satisfy, for any $R>0$, a weak large deviation principle 
on some appropriate set of sufficiently regular functions $[-R,R]^d\to(0,\infty)$ with speed
$\alpha_n^d b_n^q$ and rate function
$$
\Phi_{D,q}(\varphi)=D\int_{[-R,R]^d} \varphi^q(x)\,\d x,
$$
as the following heuristic calculation suggests:
\begin{equation}\label{LDP2heur}
\begin{aligned}
\P(\overline Y_n\approx \varphi\mbox{ on }[-R,R]^d)&\approx 
\P\Bigl(Y(z)>b_n \varphi\bigl({\textstyle{\frac z{\alpha_n}}}\bigr)
\mbox{ for }z\in [-R\alpha_n,R\alpha_n]^d\cap\Z^d\Bigr)\\
&\approx \prod_{z\in [-R\alpha_n,R\alpha_n]^d\cap\Z^d}
\exp\Bigl\{-D \bigl[b_n \varphi\bigl({\textstyle{\frac z{\alpha_n}}}\bigr)\bigr]^q\Bigr\}\\
&\approx \exp\Bigl\{ -D \alpha_n^d b_n^q\int_{[-R,R]^d} \varphi^q(x)\,\d x\Bigr\}.
\end{aligned}
\end{equation}
Note that the speeds of the two large deviation principles in 
\eqref{LDP1heur} and \eqref{LDP2heur} are equal because of \eqref{alphachoice}.
Using the two large deviation principles and \eqref{Znscale}, we see that
$$
\P({\smfrac 1n} Z_n>b_n)\approx \exp\Bigl\{-\frac n{\alpha_n^2} \widetilde K_{D,q}\Bigr\},
$$
where 
\begin{equation}\label{Ktildedef}
\widetilde K_{D,q}=\inf\{\Ical(\psi^2)+D\|\varphi\|_q^q\colon \psi^2\in\Fcal,
\varphi\in \Ccal_+(\R^d),\langle\psi^2,\varphi\rangle=1\Bigr\}.
\end{equation}
It is an elementary task to evaluate the infimum on $\varphi$ and to check that indeed
$K_{D,q}=\widetilde K_{D,q}$. This ends the heuristic explanation of Theorem~\ref{inter}.

The situation in the large deviation case, Theorem~\ref{lin}, 
is similar, when we put $b_n=1$. See \cite{AC02} for a heuristic argument in this case. 

We distinguish the two cases of  very large deviations (V)
and large deviations (L). The choices of $b_n$ and $\alpha_n$ in the respective cases are
the following.
\begin{equation}\label{bnalphachoice}
\begin{array}{lllrcl}
\mbox{case (V):}&\quad \mbox{Hypothesis of Theorem \ref{inter}}, &\quad 1\ll b_n\ll n^{\frac 1q},\qquad &\alpha_n&=&n^{\frac 1{d+2}}b_n^{-\frac q{d+2}},\\
\mbox{case (L):}&\quad \mbox{Hypothesis of Theorem \ref{lin}}, &\quad b_n=1,\qquad &\alpha_n&=&n^{\frac 1{d+2}}.
\end{array}
\end{equation}

\subsection{Small deviations for Gaussian sceneries}\label{Sec-Gauss}

\noindent Theorems \ref{inter} and \ref{lin}  do not handle sequences $(b_n)_n$ satisfying $a_n^{\ssup{0}}\ll b_n\ll 1$, where we recall from \eqref{weakconv} that $a_n^{\ssup{0}}$ is the scale of the convergence in distribution. In this regime, we present a partial result for Gaussian sceneries and simple random walk in $d=2$. This result is based on a deep result by Brydges and Slade \cite{BSl} about exponential moments of the renormalized self-intersection local time of simple random walk.

\begin{lemma}[Small deviations for Gaussian sceneries]\label{small} Assume that $Y(0)$ is a standard Gaussian random variable and that $(S_n)_n$ is the simple random walk, and assume that $d= 2$. Let 
${n}^{-1/2}(\log n)^{1/2} = a_n^{\ssup{0}} \ll b_n\ll a_n^{\ssup{1}} \equiv n^{-1/2} \log n$, then
\begin{equation}\label{smallasy}
\lim_{n\to\infty}\frac {\log n}{b_n^2 n}\log\P({\smfrac 1n }Z_n>b_n)= -\frac \pi 4.
\end{equation}
\end{lemma}

\begin{Proof}{Proof}
As we mentioned in Section~\ref{sec-model}, the distribution of the random walk in random scenery, $Z_n$, is easily identified in terms of the walk's {self-intersection local time} $\Lambda_n$ defined in \eqref{Lambdadef}. More precisely, the conditional distribution of $Z_n$ given the walk $S$ is $\Ncal\times \sqrt{\Lambda_n}$, where $\Ncal$ is a standard normal variable, independent of the walk. The typical behavior of the self-intersection local time is as follows \cite{BSl}
\begin{equation}\label{typical}
\E\bigl[\Lambda_n\bigr]\sim\frac 2 \pi\bigl(na_n^{\ssup{0}}\bigr)^2=\frac 2 \pi n\log n, \qquad n\to\infty.
\end{equation}
We prove now the upper bound in {\eqref{smallasy}}. Recall that $d=2$ and introduce the centered and normalized self-intersection local
time,
$$
    \gamma_n=\frac 1n \Bigl(\Lambda_n-\E\bigl[\Lambda_n\bigr]\Bigr).
$$

\noindent  Use Chebyshev's
inequality to obtain, for any $\theta>0$ and any $n\in\N$,
\begin{equation}\label{Chebychev}
    \P(\smfrac 1n Z_n>b_n) \le \E\bigl[e^{\theta Z_n}\bigr] e^{-\theta b_n n}.
\end{equation}

Using the above characterization of the distribution of $Z_n$, we see that 
\begin{equation}\label{Znident}
\E\bigl[e^{\theta Z_n}\bigr]=\E\bigl[\E\bigl[e^{\theta Z_n}\,\big|\,S\bigr]\bigr]=\E\bigl[\E\bigl[\exp\bigl\{\theta \Ncal \sqrt{\Lambda_n}\bigr\}\,\big|\, S\bigr]\bigr]=\E\bigl[e^{\frac 12\theta^2\Lambda_n}\bigr]=\E\bigl[e^{\frac 12\theta^2 n\gamma_n}\bigr]e^{\frac 12\theta^2 \E[\Lambda_n]}.
\end{equation} 
According to Theorem~1.2 in \cite{BSl}, $\lim_{n\to\infty}
\E[e^{c\gamma_n}]$ exists and is finite for any $c<c_0$, where $c_0>0$ is some
positive constant. Now
pick $\theta=\theta_n= \pi\, b_n/(2\log n)$. Note that
$\theta_n^2\, n\to 0$ because of $b_n\ll n^{-1/2}\log n$, and
therefore the first factor on the right hand side of \eqref{Znident} is bounded,
according to the above mentioned result of Brydges and Slade. Use \eqref{typical} on the right  hand side of \eqref{Znident} and substitute in \eqref{Chebychev} to obtain
$$
    \log\P({\smfrac 1n} Z_n>b_n)\le - (1+o(1)) {\frac \pi 4} \,
    \frac{b_n^2 n}{\log n}.
$$

\noindent This is the upper bound in (\ref{smallasy}).

Now we prove the lower bound in \eqref{smallasy}. Using the above characterization of the distribution of $Z_n$, we obtain, for any $\theta>0$,
\begin{equation}
    \label{loweresti1}
    \P(\smfrac 1n Z_n>b_n)\ge \P(\Ncal>\theta)\, \P\left( \Lambda_n>\frac{n^2b_n^2}{\theta^2}\right) .
\end{equation}
Fix an arbitrary $c\in (0, {\frac 2 \pi})$. We apply
\eqref{loweresti1} to
$\theta=b_n(\frac{n}{c\log n})^{1/2}$ and obtain
$$
    \log \P(\smfrac 1n Z_n>b_n)\ge -\frac 1{2}\,b_n^2 \frac{n}{c \log n}(1+o(1))
    +\log \P\left( \Lambda_n >c n\log n
    \right),\qquad n\to\infty.
$$
By the Paley--Zygmund inequality (Kahane \cite{K85} p.~8) stating
that $\P (X> r \E [X]) \ge (1-r)^2 \E [X]^2/\E [X^2]$ for all
$r\in (0,1)$ and all square-integrable random variables $X$, we
obtain that
$$
\P\left( \Lambda_n>c
n\log n \right)\geq \bigl(1-\bigl({\smfrac{c\pi}2}\bigr)^2\bigr)\,\frac{\E[\Lambda_n]^2}{\E[\Lambda_n^2]}.
$$
Recall from (\ref{typical}) that
$\E[\Lambda_n] \sim {\frac 2 \pi} n\log n$ as $n\to
\infty$. On the other hand, Bolthausen \cite{B89} proved that
Var$[\Lambda_n] = \Ocal (n^2)$. Therefore, $\E[\Lambda_n^2]\sim \E[\Lambda_n]^2$, and, consequently,
$$
\liminf_{n\to \infty}\, \P\left( \Lambda_n>c
n\log n \right) >0.
$$

\noindent Therefore,
$$
\liminf_{n\to \infty}\, \frac {\log n} {b_n^2n} \log \P(\smfrac 1n Z_n>b_n) \ge
-{\frac 1 {2c}}.
$$

\noindent Letting $c\uparrow{\frac 2 \pi}$, this yields the lower bound in (\ref{smallasy}).
\end{Proof}
\qed

\section{Variational formulas}\label{sec-varform}

In this section we prove Proposition~\ref{constpos} and Lemma~\ref{constasy}. In Section~\ref{subs:positivity} we prove a necessary and sufficient criterion for positivity of the constant $\chi_{d,p}$ defined in \eqref{chiident}. The relation to the Gagliardo-Nirenberg constant is discussed in Section~\ref{bem-GagNir}, and the relation to the constant $K_{D,q}$ defined in \eqref{kdq} is proved in Section~\ref{sec-Kchirel}, where we also finish the proof of Proposition~\ref{constpos}. Finally, Lemma~\ref{constasy} is proved in Section~\ref{Sec-ScalLim}.

\subsection{Positivity of $\boldsymbol{\chi_{d,p}}$}
\label{subs:positivity}

\begin{lemma}\label{chipos}
The constant $\chi_{d,p}$ is positive if and only if $d\leq\frac{2p}
{p-1}$.
\end{lemma}

\begin{Proof}{Proof} Certainly, it suffices to do the proof only in 
the case where $\frac 12\Gamma$ is the identity matrix. 

See \cite[Sect.~2]{C03} for an alternate proof of the positivity of $\chi_{d,p}$ in the subcritical dimensions, $d<\frac{2p}{p-1}$, using the relation to the Gagliardo-Nirenberg constant, which we explain in Section~\ref{bem-GagNir}.

Let us recall standard Sobolev inequalities (see \cite[Theorems~8.3, 
8.5]{LL97}). There are positive constants
$S_d$ for $d\geq 3$ and $S_{2,r}$ for $r>2$ such that 
\begin{equation}\label{Sobolev}
\begin{array}{rcll}
S_d\|\psi \|_{2d/(d-2)}^2&\leq &\|\nabla \psi \|_2^2,\qquad&\mbox{for }
d\geq 3, \psi \in D^1(\R^d)\cap L^2(\R^d),\\
S_{2,r}\|\psi \|_{r}^2&\leq& \|\nabla \psi \|_2^2+\|\psi \|_2^2,
\qquad&\mbox{for }d=2, \psi \in H^1(\R^d),\; r>2.
\end{array}
\end{equation} 
Here $D^1(\R^d)$ denotes the set of locally integrable functions $\R^d\to\R$
which vanish at infinity and possess a distributional derivative in $L^2(\R^d)$.

Let us first do the proof for the case $3\le d\leq\frac{2p}
{p-1}$.
For any $\psi\in H^1(\R^d)$ that satisfies 
$\|\psi\|_2=1=\|\psi\|_{2p}$, 
we may use the above Sobolev inequality and obtain that $\|\nabla 
\psi\|_2^2
\geq \const\|\psi\|_{2d/(d-2)}^2$. We now rewrite 
$$
\int_{\R^d} \psi^{\frac {2d}{d-2}}(t)\,\d t=\int_{\R^d}\bigl(\psi^{2p-2}(t)\bigr)^{\frac
2{(d-2)(p-1)}}\,\psi^2(t)\,\d  t.
$$
Recall that $\psi^2$ is a probability density. Therefore, an application 
of 
Jensen's inequality to the convex map $x\mapsto x^{2/[(d-2)(p-1)]}$ yields 
that 
$\|\psi\|_{2d/(d-2)}$ satisfies a lower bound in terms of a power of $\|\psi\|_{2p}$, 
which is equal to one. Hence, on the set of those $\psi\in H^1(\R^d)$ 
that satisfy $\|\psi\|_2=1=\|\psi\|_{2p}$, the map $\psi\mapsto
\|\nabla \psi\|_2^2$ is bounded away from zero. Now compare to 
\eqref{chiident} 
to see that this implies the assertion in the case $3\le d \le \frac{2p} {p-1}$.

Now we turn to $d=2$ with $p>1$ arbitrary. By a scaling $\psi_\beta= \beta^{\frac 
d2}\psi(\cdot\,\beta)$, we can
find, for any $\delta>0$, a $c(\delta)>0$ such that
\begin{equation}\label{chid=2}
\chi_{2,p}=c(\delta)\inf\Bigl\{\|\nabla \psi\|_2^2\colon \psi\in  H^1(\R^d)
,\|\psi\|_2=1,\|\psi\|_{2p}=\delta\Bigr\}.
\end{equation}
Now we choose $\delta$ such that $2\delta^{-2}=S_{2,2p}$, the Sobolev 
constant in \eqref{Sobolev}
for $d=2$ and $r=2p$. Then we have, for any $\psi$ in the set on the right 
hand side of \eqref{chid=2},
$$
2=\frac 2{\delta^2}\|\psi\|_{2p}^2=S_{2,2p}\|\psi\|_{2p}^2\leq \|\nabla \psi 
\|_2^2+\|\psi \|_2^2=\|\nabla \psi \|_2^2+1,
$$
and hence it follows that $\chi_{2,p}\geq c(\delta)>0$.

Now we show that $\chi_{2d,p}\leq 2\chi_{d,p}$ for any $d\in\N$ and 
$p\in (0,\infty)$. This simply 
follows from the observation that, for any $\psi\in H^1(\R^d)$, the 
function $\psi\otimes \psi\in H^1(\R^{2d})$
satisfies
$$
\| \nabla(\psi\otimes \psi)\|_2^2=2 \|\nabla\psi\|_2^2.
$$
Using this, the estimate $\chi_{2d,p}\leq 2\chi_{d,p}$ easily follows, 
since $\|\psi\otimes\psi \|_2=\|\psi\|^2_2$ and $\|\psi\otimes \psi\|_{2p}=\|\psi\|^2_{2p}$. 
In particular, this shows that $\chi_{1,p}>0$ for any $p>1$.

It remains to show that $\chi_{d,p}=0$ for $d>\frac{2p}{p-1}$. It is sufficient to 
construct a sequence of sufficiently regular functions $\psi_n\colon\R^d\to[0,\infty)$ 
such that $\|\psi_n\|_2$ and $\|\psi_n\|_{2p}$ both converge towards some positive numbers,
but $\|\nabla\psi_n\|_2$ vanishes as $n\to\infty$. In order to do this, pick some 
rotationally invariant function $\psi^2=f\circ|\cdot|\in\Fcal$ whose radial part 
$f\colon(0,\infty)\to(0,\infty)$ satisfies
$$
f(r)=D\times\begin{cases}r^{-\gamma} &\mbox{if 
}r\in(0,1),\\
1&\mbox{if }r\in[1,A],\\
A^{2d} r^{-2d}&\mbox{if }r>A,
\end{cases}
$$
where $A,D,\gamma>0$ are constants to be determined. Let $\omega_d$ denote the surface of the unit ball in $\R^d$. 
The following statements can be easily verified by some tedious but elementary calculations: 
\begin{eqnarray}
\gamma<d&\Longrightarrow&\|\psi\|_2^2=\frac{\omega_d}d D\Bigl[2A^d+\frac \gamma{d-\gamma}\Bigr]<\infty,\\
\gamma<\frac dp&\Longrightarrow&\|\psi\|_{2p}^{2p}={\omega_d} D^p\frac pd \Bigl[\frac\gamma{d-p\gamma}+A^d\frac 2{2p-1}\Bigr]<\infty,\\
\gamma<d-2&\Longrightarrow&\|\nabla \psi\|_2^2=\frac 14 \omega_d D\Bigl[\frac{\gamma^2}{d-\gamma-2}+A^{d-2}\frac{4d^2}{2+d}\Bigr]<\infty.
\end{eqnarray}
Since $p>1$ and $\frac dp<d-2$, we only have to assume that $\gamma<\frac dp$. Now we pick sequences $D_n$, $A_n$ and $\gamma_n$ such that all the following conditions are satisfied as $n\to\infty$:
$$
D_n\to 0,\qquad A_n\to\infty,\qquad \gamma_n\uparrow \frac dp,\qquad D_nA_n^d\to 1,\qquad\frac{D^p_n}{d-p\gamma_n}\to 1.
$$
Let $\psi_n$ be defined as the $\psi$ above with these parameters. Then we have, as $n\to\infty$,
$$
\|\psi_n\|_2^2\to 2\frac{\omega_d}d,\qquad \|\psi_n\|_{2p}^{2p}\to \omega_d,\qquad \|\nabla \psi_n\|_2^2\to 0.
$$
This ends the proof.
\end{Proof}
\qed

\subsection{Relation to the Gagliardo-Nirenberg constant}\label{bem-GagNir}

\noindent Actually, for dimensions $d\geq 2$ in the special case that $\frac 12 \Gamma$ is the identity matrix, the constant $\chi_{d,p}$ in \eqref{chiident} can be identified in terms of the {\it Gagliardo-Nirenberg\/} constant, $\kappa_{d,p}$, as follows. Assume that $d\geq 2$ and $1<p<\frac d{d-2}$. Then $\kappa_{d,p}$ is defined as the smallest constant $C$ in the {\it Gagliardo-Nirenberg inequality\/} 
\begin{equation}\label{GagNir}
\|\psi\|_{2p}\leq C \|\nabla\psi\|_2^{\frac{d(p-1)}{2p}}\|\psi\|_2^{1-\frac{d(p-1)}{2p}},\qquad \psi\in H^1(\R^d).
\end{equation}
This inequality received a lot of interest from physicists and analysts, and it has deep connections to Nash's inequality and logarithmic Sobolev inequalities. Furthermore, it also plays an important role in recent work of Chen \cite{C03} on self-intersections of random walks. See \cite[Sect.~2]{C03} for more on the Gagliardo-Nirenberg inequality.

It is clear that
\begin{equation}\label{kappadef}
\kappa_{d,p}=\sup_{\psi\in H^1(\R^d),\psi\not=0}\frac{\|\psi\|_{2p}}{\|\nabla\psi\|_2^{\frac{d(p-1)}{2p}}\|\psi\|_2^{1-\frac{d(p-1)}{2p}}}
=\Bigl(\inf_{\psi\in H^1(\R^d)\colon\|\psi\|_2=1}\|\psi\|_{2p}^{-\frac {4q}d}\|\nabla\psi\|_2^2\Bigr)^{-\frac d{4q}}.
\end{equation}
Clearly, the term over which the infimum is taken remains unchanged if $\psi$ is replaced by $\psi_\beta(\cdot)=\beta^{\frac d2}\psi(\cdot\,\beta)$ for any $\beta>0$. Hence, we can freely add the condition $\|\psi\|_{2p}=1$ and obtain that $\kappa_{d,p}=\chi_{d,p}^{-\frac d{4q}}$. In particular, the variational formulas for $\kappa_{d,p}$ in \eqref{kappadef} and for $\chi_{d,p}$ in \eqref{chiident} have the same maximizer(s) respectively minimizer(s). It is known that \eqref{kappadef} does possess a maximizer, and this is an infinitely smooth, positive and rotationally invariant function (see \cite{We83}). Uniqueness of the minimizer holds in $d\in\{2,3,4\}$ for any $p\in(1,\frac d{d-2})$, and in $d\in\{5,6,7\}$ for any $p\in(1,\frac 8d)$, see \cite{MS81}.

\subsection{Relation between $\boldsymbol{K_{D,q}}$ and $\boldsymbol{\chi_{d,p}}$ (Proposition~\ref{constpos})}\label{sec-Kchirel}

Now we prove the remaining assertions of Proposition~\ref{constpos}.

(i) The relation \eqref{Kdef} is proved by an elementary scaling argument and optimization.
Indeed, replace $\psi$ by $\psi_\beta(\cdot)=\beta^{d/2}\psi(\cdot\,\beta)$ in \eqref{kdq} and 
optimize explicitly on $\beta>0$. Afterwards the additional constraint $\|\psi\|_{2p}=1$
may freely be added. {From} \eqref{Kdef} and Lemma~\ref{chipos} the last assertion follows.

(ii) We only show the positivity of $K_H(u)$ for $d\geq 3$ and 
${\overline p}<\frac d{d-2}$; the argument for $d\leq 2$ and any ${\overline p}>1$ is the same. 

Since we assumed that $\E[Y(0)]=0$, we may pick some 
$\delta>0$ such that $H(t)\leq u t/2$ for 
$t\in[0,\delta]$. 
Pick $\eps>0$ such that ${\overline p}+\eps<
\frac d{d-2}$, then there is $c(\delta,\eps)>0$ 
depending on 
$\delta,\eps$ and $H$ only, such that $H(t)\leq c(\delta,\eps) t^{{\overline p}+\eps}$ for any
$t\in[\delta,
\infty)$. Then $H(t) \le \frac u2 \, t + c(\delta,\eps) t^{{\overline p}+\eps}$
for any $t\ge 0$, which implies that, for any $\psi\in H^1(\R^d)$ satisfying $\|\psi\|_2=1$,
$$
\begin{aligned}
\Phi_H(\psi^2,u)&\geq \sup_{\gamma>0}\Bigl[\gamma u-\int 
\frac{u\gamma}2\, \psi^2(x)\,\d x-\int c(\delta,\eps) (\gamma 
\psi^2(x))^{{\overline p}+\eps}\,\d x\Bigr]\\
&= \sup_{\gamma>0}\Bigl[ \, \frac u 2 \,
\gamma -c(\delta,\eps)\gamma^{{\overline p}+\eps}
\|\psi^2\|_{{\overline p}+\eps}^{{\overline p}+\eps}\Bigr].
\end{aligned}
$$
Now carry out the optimization over $\gamma$ to see that
$$
\Phi_H(\psi^2,u)\geq C \|\psi^2\|_{{\overline p}+\eps}^{-q_\eps},\qquad \mbox{where }\frac1{{\overline p}+\eps}+\frac 1 {q_\eps}=1,
$$
and $C>0$ depends on $u$, ${\overline p}+\eps$ and $c(\delta,\eps)$ only.
Hence, $K_H(u)\geq K_{C,q_\eps}$. Since $d \le 2q_\eps$, this is positive by assertion (i).

Now we show that $K_H(u)=0$ for any $u>0$ if ${\underline p}\equiv \liminf_{t\to\infty}
\frac{\log H(t)}{\log t}>\frac d{d-2}$. First we do this for a random variable $\widetilde Y(0)$ under the  
assumption that $\E[\widetilde Y(0)]=1$. Pick $\eps>0$ such that ${\underline p}-\eps>\frac d{d-2}$. Since $H'(0)=1$,
there is $C>0$ such that $H(t)\geq C t^{{\underline p}-\eps}$ for any $t\geq 0$. Hence, the above argument applies
and shows that $K_H(u)\leq K_{D,q_\eps}$ for some $D>0$, where $q_\eps$ is determined by
$\frac 1{{\underline p}-\eps}+\frac 1{q_\eps}=1$. Since ${\underline p}-\eps>\frac d{d-2}$, the condition $d\leq \frac{2({\underline p}-\eps)}{{\underline p}-\eps-1}$ is violated. 
Again assertion (i) implies that $K_H(u)=0$.

Let now $Y(0)$ have expectation $0$, then $\widetilde Y(0)=Y(0)+1$ has expectation 1. If $\widetilde H$ denotes the cumulant generating function of $\widetilde Y(0)$, then we have, according to the above, $K_{\widetilde H}(u)=0$ for any $u>0$. Since $K_{\widetilde H}(u)$ is well-defined, non-negative and non-decreasing for all $u\in\R$, we also have $K_{\widetilde H}(u)=0$ for any $u\in \R$. Obviously, $\widetilde H(t)=H(t)+t$ and $K_{\widetilde H}(u)=K_H(u-1)$ for any $u\in\R$, and this implies the statement.

\subsection{Scaling relations (Lemma~\ref{constasy})}\label{Sec-ScalLim}

In this section, we prove Lemma~\ref{constasy}. 

(i) Fix $\eps>0$, then there is some $C>0$ such that
$$
-Ct+(\widetilde D-\eps)t^p\leq H(t)\leq Ct+(\widetilde D+\eps)t^p,\qquad t\geq 0.
$$
Using this in the definition of $\Phi_H(\psi^2,u)$, we obtain, for any $\psi\in H^1(\R^d)$,
$$
\sup_{\gamma>0}\Bigl\{\gamma (u-C)-\gamma^p(\widetilde D+\eps)\|\psi^2\|_p^p\Bigr\}
\leq \Phi_H(\psi^2,u)\leq \sup_{\gamma>0}\Bigl\{\gamma (u+C)-\gamma^p
(\widetilde D-\eps)\|\psi^2\|_p^p\Bigr\}.
$$
The suprema may easily be evaluated, and we obtain, for some $\eta_1, 
\eta_2>0$, which vanish as $\eps\downarrow 0$,
$$
(D-\eta_1)\|\psi^2\|_p^{-q} (u-C)^q\leq \Phi_H(\psi^2,u)\leq (D+\eta_2)\|\psi^2\|_p^{-q} (u+C)^q.
$$
Using this in the definition of $K_H(u)$ in \eqref{K3def}, we obtain
$$
K_{(D-\eta_1)(u-C)^q,q}\leq K_H(u)\leq K_{(D+\eta_2)(u+C)^q,q}.
$$
Now use Proposition~\ref{constpos}(i), in particular \eqref{Kdef}, and 
use that $\eta_1,\eta_2\to 0$ as $\eps\downarrow0$.

(ii) Substituting $\psi(\cdot)=u^{d/(d+2)}\psi_0(\cdot\, u^{2/(d+2)})$ 
and $\gamma= u^{(2-d)/(2+d)}\gamma_0$ yields that 
\begin{equation}\label{KHscal}
u^{-\frac 4{d+2}}K_H(u)= \inf_{\|\psi_0\|_2=1}\Bigl\{\frac 12 
\|\Gamma^{\frac 12}\nabla \psi_0\|_2^2 +\sup_{\gamma_0>0}\Bigl( 
\gamma_0-\int u^{-2}H\bigl(u\gamma_0 \psi_0^2(x)\bigr)\,\d x\Bigr)\Bigr\}.
\end{equation}
It remains to show that the limit superior of the right 
hand side as $u\downarrow 0$ is not larger than $K_{\frac 12,2}$. 
This is shown as follows. Let $\psi_*\in H^1(\R^d)$ be an 
$L^2$-normalized bounded minimizer in the variational formula in 
\eqref{kdq} for $D=\frac 12$ and $q=2$. Its existence is proven in the same way as in
\cite{We83}, where the case $\Gamma={\rm Id}$ was considered. 
Hence we have $\lim_{u\downarrow 0}\int u^{-2}H\bigl(u\gamma_0 \psi_*^2(x)\bigr)\,\d x
=\frac 12 \gamma_0^2 \|\psi_*^2\|_2^2$, uniformly in $\gamma_0$ on compacts of $[0,\infty)$. 
Hence, the supremum on the right hand side of \eqref{KHscal} converges towards 
$\sup_{\gamma_0>0}(\gamma_0-\frac 12 \gamma_0^2 \|\psi_*^2\|_2^2)=\frac 12 \|\psi_*^2\|_2^{-2}$. 
Replacing on the right hand side of \eqref{KHscal} the infimum on $\psi_0$ by $\psi_*$, 
we arrive at $\limsup_{u\downarrow 0}u^{-4/(d+2)} K_H(u)\leq K_{\frac 12,2}$, which is 
\eqref{Kasy3}.

\section{Proof of Theorems~\ref{inter} and \ref{lin}: Preparations}\label{sec-prep}

\noindent In this section we prepare for the proofs of our main
results, Theorems~\ref{inter} and \ref{lin}. Our proofs follow the 
strategy of the proof of \cite[Theorem~2.2]{AC02}. That is, the proofs
of the lower bounds essentially follow the outline described in Section~\ref{sec-heur},
and the proofs of the upper bounds use an exponential Chebyshev inequality with a random 
parameter. However, due to the unboundedness of the scenery in our case, we face a 
serious additional difficulty, which we will overcome using a recently developed technique.

As we have already indicated in Section~\ref{sec-heur}, our main tools are large deviation 
principles for the walker's local times and for the scenery. These principles are presented
in Sections~\ref{sec-LDP} and \ref{sec-LDPsc}, respectively.
However, for the application of these two principles, there are three main technical obstacles: 
\begin{enumerate}
\item[(1)] the principles hold only on compact subsets of the space, 
\item[(2)] the scaled scenery must be smoothed,  
\item[(3)] the scaled scenery must be cut down to bounded size. 
\end{enumerate}

The first obstacle will be handled later by making a connection to the 
periodized version of the random walk, which is a standard recipe.
Hence, it will be necessary to approximate the variational formulas appearing in our main results
by finite-space versions, and this is carried out in Section~\ref{sec-appr}.
The necessity of the smoothing arises from the fact that the map $(\psi^2,\varphi)\mapsto
\langle \psi^2,\varphi\rangle$ is not continuous in the product of the topologies on which 
the large deviation principles are based. This was already pointed out in \cite{AC02}. The remedy is a smoothing
procedure which was introduced in \cite{AC02} and will be adapted in Section~\ref{sec-smooth} below.
However, this procedure only works for {\it uniformly bounded\/} sceneries, and this 
explains the necessity of a cutting argument for the scenery. This obstacle was not present
in \cite{AC02} and is the main technical challenge in the present paper, see Section~\ref{sec-cut}.

\subsection{Large deviations for the local times}\label{sec-LDP}

\noindent In this section, we formulate one of our main tools: large deviation 
principles for the normalized and scaled local times.
These principles are essentially standard and well-known, however, some of the principles
we use do not seem to have been proven in the literature, and therefore we shall provide
a proof for them in the appendix.

For the convenience of the reader, we recall the notion of a large deviation
principle. A sequence $(X_n)_{n\in\N}$ of random variables (or their distributions),
taking values in a topological space $\Xcal$,
satisfy a {\it large deviation principle\/} with {\it speed\/} $(\gamma_n)_{n\in\N}$ and 
{\it rate function\/} $\Ical\colon \Xcal\to[0,\infty]$, if the following two statements hold:
\begin{eqnarray}
\limsup_{n\to\infty}\frac1{\gamma_n}\log P(X_n\in 
F)&\leq& - \inf_{F} \Ical,\qquad F\subset \Xcal\mbox{ 
closed},\label{upperbound}\\
\liminf_{n\to\infty}\frac1{\gamma_n}\log P(X_n\in 
O)&\geq& -\inf_{O} \Ical,\qquad O\subset \Xcal\mbox{ 
open}.\label{lowerbound}
\end{eqnarray}
This definition equally applies if the measure $P$ has not full mass, but happens to be 
a subprobability measure only.

We shall need large deviation principles for a rescaled version of the local times of 
our random walk. More precisely, we shall 
need two slightly different principles: 
one on never leaving a given cube in $\Z^d$ and $\R^d$, respectively, and another one for the 
periodized version of the walk on that cube. 
We recall
that we have listed our assumptions on the random walk at the beginning of Section~\ref{results}.
For $R>0$, we denote by $B_R=[-R,R]^d\cap\Z^d$ the centered box in
$\Z^d$ with radius $R$. By $S^{\smallsup{R}}=(S^{\smallsup{R}}_0,S^{\smallsup{R}}_1,\dots)$
we denote the random walk on the torus $B_R$, i.e., the walk on $B_R$ (with the 
opposite sides identified with each other) having transition kernel
\begin{equation}\label{kernelR}
p^{\smallsup{R}}(z,\widetilde z)=\sum_{k\in\Z^d}p(z,\widetilde z+2k\lfloor 
R\rfloor),\qquad z,\widetilde z\in B_R,
\end{equation}
where $p(\cdot,\cdot)$ denotes the transition kernel of $S$. Note 
that $p^{\smallsup{R}}$ is symmetric since $p$ is. 
The local times of $S^{\smallsup{R}}$ are denoted by
\begin{equation}\label{perloctim}
\ell_n^{\ssup{R}}(z)=\sum_{k\in\Z^d}\ell_n(z+2k\lfloor R\rfloor),\qquad z\in B_R.
\end{equation}

We consider rescaled versions of $\frac 1n\ell_n$ and
$\frac 1n \ell_n^{\smallsup{R}}$. Recall the normalized and rescaled  version 
$L_n$ of the local times $\ell_n$ defined in \eqref{Lndef}. 
By $\Fcal_R$ we denote the subset of those functions in $\Fcal$ whose 
support lies in $Q_R=[-R,R]^d$. 
Note that 
\begin{equation}\label{supps}
\supp(L_n)\subset Q_R\qquad\Longleftrightarrow 
\qquad\supp(\ell_n)\subset B_{R\alpha_n}.
\end{equation}

Denote the scaled version of the torus-version of the local times, $\frac 1n \ell_n^{\smallsup {R\alpha_n}}$, by
$L_n^{\smallsup{R}}\colon Q_R\to[0,\infty)$. Then 
$L_n^{\smallsup{R}}$ is a random element of the set $\Fcal^{\smallsup{R}}$ of probability densities on the 
torus $Q_R=[-R,R]^d$, whose opposite sides are identified with each other.
We define a rate function 
$\Ical^{\smallsup{R}}\colon\Fcal^{\smallsup{R}}\to[0,\infty]$ by
\begin{equation}\label{IRdef}
\Ical^{\smallsup{R}}(\psi^2)=\frac 1{2}\int_{Q_R}\bigl|\Gamma^{\frac 12}\nabla_R
\psi(x)\big|^2\,\d x,
\end{equation}
if $\psi$ has an extension to an element of $H^1(\R^d)$, 
and $\Ical^{\smallsup{R}}(\psi^2)=\infty$ otherwise. Here $\nabla_R$ 
denotes the gradient on the torus $Q_R$, i.e., with periodic boundary condition.

The topology used on the sets $\Fcal_R$
and on $\Fcal^{\smallsup{R}}$ are the weak topologies induced by the 
test integrals against 
the continuous bounded functions on
$Q_R$. If we identify any element of $\Fcal_R$ resp.~of
$\Fcal^{\smallsup{R}}$
with a probability measure, then this topology is just the usual weak 
topology
on the set of probability measures on $Q_R$. In this
case, we extend the respective rate functions trivially by $\infty$ to 
the set of measures not having a density.

\begin{lemma}[Large deviation principles for $L_n$]\label{LDP} Fix $R>0$. 
Assume that $\alpha_n \to \infty$ and 
$$
\alpha_n^d\ll \begin{cases} \sqrt n&\mbox{if }d=1,\\
\frac n{\log n}&\mbox{if }d=2,\\
n&\mbox{if }d\geq 3,
\end{cases}
$$
 as $n\to\infty$.  Then 
the following two facts hold true.
\begin{enumerate}
\item[(i)] The distributions of $L_n$ under $\P(\,\cdot\,\cap\{\supp(L_n)\subset Q_R\})$
satisfy a  large deviation principle on $\Fcal_R$
with speed $n\alpha_n^{-2}$ and rate function $\Ical_R$, the 
restriction of $\Ical$ defined in 
\eqref{Idef} to $\Fcal_R$.

\item[(ii)] The distributions of $L_n^{\smallsup{R}}$ under $\P$ 
satisfy a  large deviation principle on $\Fcal^{\smallsup{R}}$ with 
 speed $n\alpha_n^{-2}$ and rate function $\Ical^{\smallsup{R}}$ 
given in \eqref{IRdef}.
\end{enumerate}
\end{lemma}

The upper bound \eqref{upperbound} of the principle in (i) for the 
special case of 
simple random walk and $\alpha_n=n^{\frac 1{d+2}}$ has been 
proven by Donsker and Varadhan \cite{DV79}, Section~3.
We have deferred the proof of Lemma~\ref{LDP} to the 
Appendix, Section~\ref{sec-proofLDP}. We feel that 
the statement and its proof are standard and should be known 
to the experts,
but we could not find a reference in the literature. 
Our proof basically follows the route of \cite{Ga77}, which has 
become standard by now. The strategy for the proof of (i) can be roughly summarized 
as follows (the proof of (ii) is analogous). We shall identify the cumulant generating function 
of $L_n$ (i.e., the logarithmic asymptotics of exponential moments 
of test integrals against continuous and bounded functions $f$) 
in terms of the Dirichlet eigenvalue of the operator 
$\frac{1}{2}\nabla\cdot\Gamma\nabla + f$. In a second step, we prove 
the large deviation principle via what is called now the 
abstract {\em G\"artner-Ellis theorem\/} and identify the 
rate function of the large deviation principle as the Legendre 
transform of the eigenvalue.

\subsection{Large deviations for the scenery}\label{sec-LDPsc}

\noindent In the proofs of the lower bounds in Theorems~\ref{inter} and \ref{lin}, we shall rely on
precise large deviation lower bounds for the scenery, tested against fixed
functions. The precise formulations are given for the respective cases here. 
Recall from \eqref{bnalphachoice} the two cases (V) and (L), which correspond to Theorems~\ref{inter} and \ref{lin}, respectively.

\noindent We begin, in case (V ) with a large deviation
principle for the rescaled scenery $\overline Y_n$ defined in \eqref{Yscaled}.

\begin{lemma}\label{LDPfieldcont} Assume the case (V) in (\ref{bnalphachoice}), and 
pick sequences $(b_n)_n$ and $(\alpha_n)_n$ as in \eqref{bnalphachoice}.
Fix $R>0$ and a continuous function
$\varphi\colon Q_R\to (0,\infty)$. Then
\begin{equation}
\liminf_{n\to\infty}\frac1{\alpha_n^d b_n^q}\log\P\bigl(\overline Y_n 
\geq \varphi\mbox{ on }Q_R\bigr)\geq-D\|\varphi\|_q^q.
\end{equation}
\end{lemma}

\begin{Proof}{Proof} Fix some small $\eps>0$. 
It is easy to see that, for sufficiently large $n\in\N$,
$$
\begin{aligned}
\P\bigl(\overline Y_n 
\geq \varphi\mbox{ on }Q_R\bigr)&=\prod_{z\in B_{R\alpha_n}}
\P\bigl(Y(z)\geq b_n \varphi(\smfrac z{\alpha_n})\bigr)\\
&\geq \exp\Bigl\{-(D-\eps)b_n^q\sum_{z\in B_{R\alpha_n}}\varphi(\smfrac z{\alpha_n})^q\Bigr\}\\
&\geq\exp\Bigl\{-(D-2\eps)\alpha_n^d b_n^q\|\varphi\|_q^q\Bigr\}.
\end{aligned}
$$
\end{Proof}
\qed

Let us now proceed with case (L). 

\begin{lemma}\label{LDPfieldL} Assume the case (L) in (\ref{bnalphachoice}) and fix $R>0$, $M>0$ and a positive continuous function $\psi^2 \colon Q_R\to (0,\infty)$. Recall that $\alpha_n=n^{\frac 1{d+2}}$. Let $\widetilde H_M$ be the conditional cumulant generating function of $Y(0)$ given that $Y(0)\geq -M$. Then, for any $u>0$,
\begin{equation}
\liminf_{n\to\infty}\frac1{\alpha_n^d}\log\P\Bigl(\int_{Q_R}\overline Y_n(x) \psi^2 (x)\,\d x\geq u\,\Big|\,Y(z)\geq -M\quad \forall z\in B_{R\alpha_n}\Bigr)\geq \Phi_{\widetilde H_M}(\psi^2 ,u;R),
\end{equation}
where 
\begin{equation}\label{PhiRdef}
\Phi_H(\psi^2,u;R)=\sup_{\gamma>0}\Bigl(\gamma u-\int_{Q_R}H(\gamma\psi^2(x))\,\d x\bigr)\Bigr)
\end{equation}
is the $Q_R$-version of $\Phi_H$ defined in \eqref{Phidef}.
\end{lemma}

\begin{Proof}{Proof} 
For any $\gamma>0$, we have
$$
\begin{aligned}
\E\Bigl[&\exp\Bigl\{\gamma\alpha_n^d\int_{Q_R}\overline Y_n(x) \psi^2 (x)\,\d x\Bigr\}\,\Big|\,Y(z)\geq -M\quad \forall z\in B_{R\alpha_n}\Bigr]\\
&=\E\Bigl[\exp\Bigl\{\gamma\sum_{z\in B_{R\alpha_n}} Y(z)\alpha_n^d\int_{z/\alpha_n+[0,1/\alpha_n]^d}\psi^2(x)\,\d x\Bigr\}\,\Big|\,Y(z)\geq -M\quad \forall z\in B_{R\alpha_n}\Bigr]\\
&=\prod_{z\in B_{R\alpha_n}}e^{(1+o(1))\widetilde H_M\big(\gamma \psi^2 ({\smfrac {z}{\alpha_n}})\big)}\\
&=\exp\Bigl\{{\alpha_n^d}\int_{Q_R} \widetilde H_M(\gamma \psi^2 (x))\,\d x\,(1+o(1))\Bigr\}.
\end{aligned}
$$
According to a variant of the G\"artner-Ellis theorem, $\int_{Q_R}\overline Y_n(x) \psi^2 (x)\,\d x$ satisfies, under conditioning on $Y(z)\geq -M$ for all $z\in B_{R\alpha_n}$, a large deviation principle on $(0,\infty)$ with speed $\alpha_n^{d}$ and rate given by the Legendre transform of the map $\gamma\mapsto \int_{Q_R} \widetilde H_M(\gamma \psi^2 (x))\,\d x$. This transform is equal to the map $u\mapsto\Phi_{\widetilde H_M}(\psi^2 ,u;R)$.
\end{Proof}
\qed

\subsection{The cutting argument}\label{sec-cut}

\noindent In this section we provide the cutting argument for the scenery in the cases (V) and  (L). 
Our method consists of a careful analysis of the $k$-th moments of the random walk in random
scenery, where $k=k_n$ is chosen in an appropriate dependence of $n$. Variants of this
method have recently been developed in the study of mutual intersections of random paths
in \cite{C03} and \cite{KM02}. 

Fix sequences $(b_n)_n$ and $(\alpha_n)_n$ as in \eqref{bnalphachoice}
and consider the scaled normalized scenery $\overline Y_n$ as defined in \eqref{Yscaled}.
Fix $M>0$.  We use the notation 
\begin{equation}\label{ydeco}
y^{\ssup{\leq M}}=(y\wedge M)\vee (-M)\qquad \mbox{and}\qquad 
y^{\ssup{> M}}=(y-M)_+,\qquad\mbox{for any }y\in \R.
\end{equation}
Later we shall estimate the scaled scenery $\overline Y_n$ by 
$\overline Y_n\leq \overline Y_n^{\ssup{\leq M}}
+\overline Y_n^{\ssup{> M}}$. Here we show how we shall handle the second term.

\begin{prop}[Scenery cutting]\label{fieldcut} Assume one of the cases (V) or (L) in (\ref{bnalphachoice}).
Then, for any $\eps>0$,
\begin{equation}
\lim_{M\to\infty}\limsup_{n\to\infty}\frac{\alpha_n^2}n\log
\P\bigl(\langle L_n,\overline Y_n^{\ssup{> M}}\rangle >\eps)=-\infty.
\end{equation}
\end{prop}

\setcounter{step}{0}

\begin{Proof}{Proof}

\bes It suffices to establish that there exists $C_M>0$ satisfying $\lim_{M\to\infty }C_M=0$ and 
\begin{equation}\label{goal}
\E\bigl[\langle \ell_n, Y^{\ssup{>Mb_n}}\rangle^k\bigr]\leq n^k b_n^k C_M^k,\qquad n\in\N,\mbox{ where }k=\frac n{\alpha_n^2}.
\end{equation}
\es

\begin{Proof}{Proof}  Use the Markov inequality to estimate, for any $\eps,M>0$ and $n,k\in\N$,
$$
\P\bigl(\langle L_n,\overline Y_n^{\ssup{> M}}\rangle>\eps)\leq \eps^{-k}\E\bigl[\langle L_n,\overline Y_n^{\ssup{>M}}\rangle^k\bigr]=\eps^{-k}(nb_n)^{-k}\E\bigl[\langle \ell_n,Y^{\ssup{>Mb_n}}\rangle^k\bigr].
$$
Now put $k=n\alpha_n^{-2}$ and observe that the estimate in \eqref{goal} for some $C_M\to 0$ as $M\to\infty$ implies Proposition~\ref{fieldcut}.
\end{Proof}\qed

Our next step is a variant of the well-known periodization technique which projects the random walk in random scenery into a fixed box. Recall from \eqref{perloctim} the local times of the periodized random walk.

\bes[Periodization]\label{step-per} For any $R,n,k\in\N$ and for any i.i.d.\ scenery $Y$ which is independent of the random walk,
\begin{equation}\label{writeoutform}
\E\bigl[\langle\ell_n,Y\rangle^k\bigr]\leq \sum_{z_1,\dots,z_k\in B_R}\E\Bigl[\prod_{i=1}^k \ell_n^{\ssup{R}}(z_i)\Bigr]\prod_{x\in B_R}\E\bigl[|Y(0)|^{\#\{i\colon z_i=x\}}\bigr].
\end{equation}
\es

\begin{Proof}{Proof} We write out
\begin{equation}\label{writeout}
\E\bigl[\langle\ell_n,Y\rangle^k\bigr]=\sum_{z_1,\dots,z_k\in B_R}\sum_{m_1,\dots,m_k\in\Z^d}\E\Bigl[\prod_{i=1}^k\ell_n(z_i+2Rm_i)\Bigr]\E\Bigl[\prod_{i=1}^k Y(z_i+2Rm_i)\Bigr].
\end{equation}
We use that the scenery is i.i.d.~and derive, with the help of Jensen's inequality, the estimate
$$
\begin{aligned}
\E\Bigl[\prod_{i=1}^k Y(z_i+2Rm_i)\Bigr]&=\prod_{x\in B_R}\prod_{y\in \Z^d}\E\Bigl[Y(y)^{\#\{i\colon z_i=x,z_i+2Rm_i=y\}}\Bigr]\\
&\leq \prod_{x\in B_R}\prod_{y\in \Z^d}\E\Bigl[|Y(0)|^{\#\{i\colon z_i=x\}}\Bigr]^{\frac{\#\{i\colon z_i=x,z_i+2Rm_i=y\}}{\#\{i\colon z_i=x\}}}\\
&=\prod_{x\in B_R}\E\Bigl[|Y(0)|^{\#\{i\colon z_i=x\}}\Bigr].
\end{aligned}
$$
Use this in \eqref{writeout} and carry out the sum over $m_1,\dots,m_k$ to finish.
\end{Proof}\qed

In the next step we estimate the term in \eqref{writeoutform} that involves the walker's local times. We denote by $S^{\ssup{R}}$ the periodized version of the random walk in $B_R$ and by $p_s^{\ssup{R}}(x,y)$ its transition probability from $x$ to $y$ in $s$ steps. By
\begin{equation}
G^{\ssup{R}}_\lambda (x, y) = \sum\limits_{s=0}^\infty e^{-\lambda s}p_s^{\ssup{R}}(x,y),
\end{equation}
we denote  the Green's function associated with the periodized walk, geometrically stopped with parameter $\lambda>0$. $\Sym_k$ denotes the set of permutations of $1,\dots,k$.

\bes\label{step-LocTim} Fix $R>0$, $\lambda>0$ and $k\in\N$. 
Then,  for any $n\in\N$, and for any $z_1,\dots,z_k\in B_{R}$,
\begin{equation}\label{locestiGreen}
\E\Bigl[\prod_{i=1}^k \ell_n^{\ssup{R}}(z_i)\Bigr]\leq e^{\lambda n}\sum_{\sigma\in\Sym_k}\prod_{i=1}^k G^{\ssup{R}}_{\lambda}\bigl(z_{\sigma(i-1)},z_{\sigma(i)}\bigr).
\end{equation}
\es

\begin{Proof}{Proof} Writing out the local times, we obtain
\begin{equation}\label{loctimesti1}
\begin{aligned}
\E\Bigl[\prod_{i=1}^k \ell_n^{\ssup{R}}(z_i)\Bigr]&\leq  \sum_{t_1,\dots,t_k=0}^n\P\Bigl(S^{\ssup{R}}_{t_i}=z_i, \, i=1,\dots,k\Big)\\
&\leq \sum_{0\leq t_1\leq t_2\leq \dots\leq t_k\leq n}\sum_{\sigma\in\Sym_k}\P\Bigl(S^{\ssup{R}}_{t_{\sigma(i)}}=z_i, \, i=1,\dots,k\Big)\\
&=\sum_{\sigma\in\Sym_k}\sum_{s_1,\dots,s_k\in\N_0}\1\Bigl\{\sum_{i=1}^ks_i\leq n\Bigr\}\prod_{i=1}^k p_{s_i}^{\ssup{R}}\bigl(z_{\sigma(i-1)},z_{\sigma(i)}\bigr),
\end{aligned}
\end{equation}
where in the last line we substituted $s_i=t_i-t_{i-1}$ and wrote $\sigma^{-1}$ instead of $\sigma$. We put $\sigma(0)=0$ and $z_0=0$. Now we estimate the indicator by
$$
\1\Bigl\{\sum_{i=1}^ks_i\leq n\Bigr\}\leq e^{\lambda n}\prod_{i=1}^k e^{-s_i \lambda}.
$$
Using this in \eqref{loctimesti1} and carrying out the sums over $s_1,\dots,s_k$, we arrive at the assertion.

\end{Proof}
\qed

In order to further estimate the Greenian term on the right of \eqref{locestiGreen}, we shall later need the following.

\bes\label{step-Green}Fix $R>0$ and $p'\in(1,\frac d{d-2})$, if $d \geq 3$, or 
$p'>1$ if $d \in \{1,2\}$. Then there is a constant $C>0$ such that, for any $n\in\N$ and any $x\in B_{R\alpha_n}$,
\begin{equation}\label{Greenesti}
\sum_{y\in B_{R\alpha_n}}G^{\ssup{R\alpha_n}}_{\alpha_n^{-2}}(x,y)^{p'}\leq C\alpha_n^{d+(2-d)p'}.
\end{equation}
\es

\begin{Proof}{Proof} For $d\leq 4$, we estimate, with the help of Jensen's inequality, and using that $p_s^{\ssup{R\alpha_n}}(x,y)$ is not bigger than one and that its sum on $y\in B_{R\alpha_n}$ equals one,
$$
\begin{aligned}
\sum_{y\in B_{R\alpha_n}}G^{\ssup{R\alpha_n}}_{\alpha_n^{-2}}(x,y)^{p'}&=\sum_{y\in B_{R\alpha_n}}\Big(\sum_{s=0}^\infty e^{-s\alpha_n^{-2}}p_s^{\ssup{R\alpha_n}}(x,y)\Big)^{p'}\\
&\leq \big(1-e^{-\alpha_n^{-2}}\big)^{p'-1}\sum_{y\in B_{R\alpha_n}}\sum_{s=0}^\infty e^{-s\alpha_n^{-2}}p_s^{\ssup{R\alpha_n}}(x,y)\\
&\leq \big(1-e^{-\alpha_n^{-2}}\big)^{p'-2}\sim\alpha_n^{4-2p'}.
\end{aligned}
$$
Now noting that $4-2p'\leq d+(2-d)p'$ for $d\leq 4$ finishes the proof of \eqref{Greenesti}.

For $d\geq 4$, we use another argument, which is based on the estimate \cite[Th.~2]{Uc98} $G(0,y)\leq C |y|^{2-d}$ for any $y\in\Z^d\setminus\{0\}$, where $G$ is the Green's function for the free (i.e., non-stopped and non-periodized) random walk, and $C>0$ is constant. Certainly, it suffices to take $x=0$. We use $C>0$ and $c >0$ to denote generic 
positive constants, not depending on $n$ or $y$, which may change their values from line to line. We estimate
\begin{equation}\label{Gsplit}
G^{\ssup{R\alpha_n}}_{\alpha_n^{-2}}(0,y)\leq G(0,y)+\sum_{m\in\Z^d\setminus\{0\}}\sum_{s\in\N_0}e^{-s\alpha_n^{-2}}p_s(0,y+2mR\alpha_n).
\end{equation}
For the first term, use the above mentioned result to see that $\sum_{y\in B_{R\alpha_n}}G(0,y)^{p'}\leq C\alpha_n^{d+(2-d)p'}$. With $\gamma>0$ a small auxiliary parameter, we split the sum on $s$ in the parts where $s\leq \gamma |m|\alpha_n$ and the remainder. Recall that the walker's steps have some exponential moments, see \eqref{emom}. Hence, we can estimate, if $\gamma$ is small enough ($\gamma<\frac R4/\log\E[e^{|S_1|}]$ suffices), for $|m|\geq 1$ and $s\leq \gamma |m|\alpha_n$, and all $y\in B_{R\alpha_n}$,
\begin{equation}\label{leqesti}
\begin{aligned}
p_s(0,y+2mR\alpha_n)&\leq \P(|S_s|\geq |y+2mR\alpha_n|)\leq \E[e^{|S_1|}]^s e^{-|y+2mR\alpha_n|}\leq \E[e^{|S_1|}]^se^{-R\alpha_n|m|}\\
&\leq e^{-c|m|\alpha_n}.
\end{aligned}
\end{equation}
This gives, for any $y\in B_{R\alpha_n}$,
\begin{equation}\label{firstsum}
\sum_{m\in\Z^d\setminus\{0\}}\sum_{s\in\N_0\colon s\leq \gamma\alpha_n |m|}e^{-s\alpha_n^{-2}}p_s(0,y+2mR\alpha_n)\leq C \sum_{m\in\Z^d\setminus\{0\}}e^{-c|m|\alpha_n}=o(\alpha_n^{2-d}).
\end{equation}
The remainder is estimated as follows. We use the local central limit theorem (see \cite[Ch.~VII, Thm.~13]{Pe75}) to deduce that there are $C>0$ and $c>0$ such that
\begin{equation}\label{LCLT}
p_s(0,x)\leq \frac C{s^{d/2}}e^{-c|x|^2/s}+C s^ {-d},\qquad s\in\N,\,  x\in\Z^d.
\end{equation}
This gives, for any $y\in B_{R\alpha_n}$,
$$
\begin{aligned}
\sum_{m\in\Z^d\setminus\{0\}}&\sum_{s\in\N_0\colon s\geq \gamma\alpha_n |m|}e^{-s\alpha_n^{-2}}p_s(0,y+2mR\alpha_n)\\
&\leq
C\sum_{ s\geq \gamma\alpha_n}e^{-s\alpha_n^{-2}}\Big[s^{-d/2}\sum_{0< |m|\leq s/(\gamma\alpha_n)}e^{-c|m|^2\alpha_n^2/s}+\Big(\frac s{\alpha_n}\Big)^d s^{-d}\Big],
\end{aligned}
$$
where we interchanged the sums on $s$ and $m$, and we also used that $|y+2mR\alpha_n|\geq |m|\alpha_n$ for $m\in\Z^d\setminus\{0\}$. Using the substitution $w=|m|\alpha_n/\sqrt s$, the sum on $m$ is estimated by
$$
\sum_{0< |m|\leq s/(\gamma\alpha_n)}e^{-c|m|^2\alpha_n^2/s}
\leq C \Big(\frac s{\alpha_n^2}\Big)^{d/2}\int_{\alpha_n/\sqrt s}^{\sqrt s/\gamma}\d w\, w^{d-1}e^{-cw^2}\leq C \Big(\frac s{\alpha_n^2}\Big)^{d/2}.
$$
Since $\sum_{ s\in\N_0}e^{-s\alpha_n^{-2}}\leq C\alpha_n^2$, this implies that 
\begin{equation}\label{geqesti}
\sum_{m\in\Z^d\setminus\{0\}}\sum_{s\in\N_0\colon s\geq \gamma\alpha_n |m|}e^{-s\alpha_n^{-2}}p_s(0,y+2mR\alpha_n)\leq C\alpha_n^{2-d}.
\end{equation}
 Use \eqref{firstsum} and \eqref{geqesti} in \eqref{Gsplit} to conclude.
\end{Proof}
\qed

The next step is a preparation for the estimate of the last term in \eqref{writeoutform}.

\bes\label{step-Yesti} Let $Y$ be a random variable that satisfies 
\begin{equation}\label{ttail}
\limsup_{r\to\infty}r^{-q}\log\P(Y>r)<0
\end{equation}
 for some $q> 1$.
\begin{enumerate}
\item[(i)]
Fix $L>0$. Then there is $C_{M,L}>0$ such that $\lim_{M\to\infty}C_{M,L}=0$ such that, for every $n\in\N$ and $M>0$ and $b_n \geq 1$,
\begin{equation}
\E\bigl[(Y-M b_n)_+^{Lb_n^q}\bigr]^{\frac 1{Lb_n^q}}\leq b_n C_{M,L}.
\end{equation}

\item[(ii)] There is a constant $C>0$ such that, for any $\mu \in\N$,
\begin{equation}
\E[Y_+^\mu ]\leq \mu^{\frac 1q \mu }C^\mu.
\end{equation}
\end{enumerate}
\es

\begin{Proof}{Proof} From our assumption on $Y$, we know that there are $C,D>0$ and $q> 1$ such that $\P(Y>s)\leq C e^{-Ds^q}$ for all $s>0$.

{\it Proof of (i).} We write $L$ instead of $Lb_n^q$ and have
\begin{equation}\label{expectY}
b_n^{-L}\E\bigl[(Y-M b_n)_+^{L}\bigr]=b_n^{-L}\int_0^\infty \P\bigl((Y-Mb_n)^L>t\bigr)\,\d t=L\int_0^\infty s^{L-1}\P(Y>(s+M)b_n)\,\d s.
\end{equation}
Now use the above estimate $\P(Y>(s+M)b_n)\leq C\exp\{-D(s+M)^q b_n^q\}$ for all $s>0$. Furthermore, use that $(s+M)^q\geq s^q+M^q$. This gives
\begin{equation}
\E\bigl[(Y-M b_n)_+^{L}\bigr] \leq b_n^L LC e^{-DM^qb_n^q}\int_0^\infty s^{L-1} e^{-D(sb_n)^q}\, \d s=LC e^{-DM^qb_n^q}\int_0^\infty s^{L-1} e^{-Ds^q}\,\d s.
\end{equation}
The change of variables $t=D s^q$ turns this into 
\begin{equation}\label{momYesti}
\E\bigl[(Y-M b_n)_+^{L}\bigr] \leq L C e^{-DM^qb_n^q}D^{-L/q}q \Gamma(L/q),
\end{equation}
where $\Gamma$ denotes the Gamma-function. 
Note that $\Gamma(x)\leq (C_1x)^x$ for some $C_1>0$ and all $x\geq 1$. Now we replace $L$ by $Lb_n^q$ and take the $(L b_n^q)$-th root to obtain 
$$
\E\bigl[(Y-M b_n)_+^{Lb_n^q}\bigr]^{\frac 1{Lb_n^q}}
\leq \widetilde C (L b_n^q)^{\frac 1{Lb_n^q}}L^{\frac 1q} e^{-M^qD/L} b_n,
$$
where $\widetilde C$ does not depend on $L$ nor on $M$ or $n$. Since $b_n\geq 1$, the assertion is proved.

{\it Proof of (ii).} From \eqref{momYesti} with $M=0$ and $L=\mu$, we have $\E[Y_+^\mu]\leq \mu C \,D^{-\mu/q} \Gamma(\mu/q)$.
Recalling that $\Gamma(x)\leq (C_1x)^x$ for some $C_1>0$ and all $x\geq 1$, we arrive at the assertion.

\end{Proof}
\qed

\bes\label{conclusion} Conclusion of the proof. 
\es
\begin{Proof}{Proof} Fix $R>0$ and let $B=B_{R\alpha_n}$ be the centered box in $\Z^d$ with radius $R\alpha_n$. 
Note that in both cases (V) and (L), (\ref{ttail}) is satisfied with $q>\frac d2$.
 Let $p$ be defined by $1=\frac 1p+\frac 1q$. Then, in both cases, $p \in (1, \frac{d}{d-2})$ if $d \geq 3$ and $p> 1$ if $d=2$. Put $k=n\alpha_n^{-2}$. Recall that $\alpha_n^{d+2}=nb_n^{-q}$. Recall that it suffices to prove \eqref{goal}. In the following, we shall use $C$ to denote a generic positive constant which depends on $R$, $q$ and $D$ only and may change its value from line to line.

Use Steps~\ref{step-per}--\ref{step-LocTim} for the scenery $Y$ replaced by $Y^{\ssup{>Mb_n}}$ and $R$ replaced by $R\alpha_n$ and with $\lambda=\alpha_n^{-2}$ to obtain
\begin{equation}\label{cuttingesti2}
\E\bigl[\langle \ell_n, Y^{\ssup{>Mb_n}}\rangle^k\bigr]\leq e^{k}\sum_{\sigma\in\Sym_k}\sum_{z_1,\dots,z_k\in B} \prod_{i=1}^k G^{\ssup{R\alpha_n}}_{\alpha_n^{-2}}\bigl(z_{\sigma(i-1)},z_{\sigma(i)}\bigr)\prod_{x\in B}\E\bigl[(Y(0)-Mb_n)_+^{\mu_x}\bigr],
\end{equation}
where we abbreviated $\mu_x=\#\{i\colon z_i=x\}$. Let us estimate the last term. We fix a parameter $L>0$ and split the product on $x\in B$ into the subproducts on $B_{\ssup{L}}=\{x\in B\colon \mu_x\leq L b_n^q\}$ and $B_{\ssup{L}}^{\rm c}=B\setminus B_{\ssup{L}}$. We estimate, with the help of Step~\ref{step-Yesti},
\begin{equation}
\begin{aligned}
\prod_{x\in B}\E\bigl[(Y(0)-Mb_n)_+^{\mu_x}\bigr]&\leq \prod_{x\in B_{\ssup{L}}}\E\bigl[(Y(0)-Mb_n)_+^{Lb_n^q}\bigr]^{\frac {\mu_x}{Lb_n^q}}\prod_{x\in B_{\ssup{L}}^{\rm c}}\E[Y(0)_+^{\mu_x}]\\
&\leq \prod_{x\in B_{\ssup{L}}}(C_{M,L}b_n)^{\mu_x}\prod_{x\in B_{\ssup{L}}^{\rm c}}(C\mu_x^{\frac 1q})^{\mu_x}.
\end{aligned}
\end{equation}
Let us abbreviate the term on the right hand side by $K(\mu)$ where $\mu=(\mu_x)_{x\in B}$. Now we pick numbers $p'>p$, $q'>1$ such that $\frac 1{p'}+\frac 1{q'}=1$ and, if $d \geq 3$,  $p'<\frac d{d-2}$, and use H\"older's inequality in \eqref{cuttingesti2} to obtain
\begin{equation}\label{cuttingesti3}
\E\bigl[\langle \ell_n, Y^{\ssup{>Mb_n}}\rangle^k\bigr]\leq e^{k}\sum_{\sigma\in\Sym_k}\Bigl(\sum_{z_1,\dots,z_k\in B} \prod_{i=1}^k G^{\ssup{R\alpha_n}}_{\alpha_n^{-2}}\bigl(z_{\sigma(i-1)},z_{\sigma(i)}\bigr)^{p'}\Bigr)^{\frac 1{p'}}\Bigl(\sum_{z_1,\dots,z_k\in B}K(\mu)^{q'}\Bigr)^{\frac 1{q'}}.
\end{equation}
Using \eqref{Greenesti} in Step~\ref{step-Green}, the term in the first brackets may be estimated by
\begin{equation}\label{Gesti}
\Bigl(\sum_{z_1,\dots,z_k\in B} \prod_{i=1}^k G^{\ssup{R\alpha_n}}_{\alpha_n^{-2}}\bigl(z_{\sigma(i-1)},z_{\sigma(i)}\bigr)^{p'}\Bigr)^{\frac 1{p'}}\leq C^k \alpha_n^{2k} \alpha_n^{-\frac 1{q'}dk}.
\end{equation}
Now we estimate the last term in \eqref{cuttingesti3}. By $A_k$ we denote the set of maps $\mu\colon B\to\N_0$ such that $\sum_{x\in B}\mu_x=k$. Observe that, for any $\mu\in A_k$, we have 
$$
\#\{(z_1,\dots,z_k)\in B^k\colon \mu_x=\#\{i\colon z_i=x\}\quad \forall x\in B\}=\frac {k!}{\prod_{x\in B}\mu_x!}.
$$
Hence,
\begin{equation}\label{betaexpression}
\sum_{z_1,\dots,z_k\in B}K(\mu)^{q'}\leq C^k k!\sum_{\mu\in A_k} \prod_{x\in B_{\ssup{L}}}C_{M,L}^{q'\mu_x}\prod_{x\in B_{\ssup{L}}}\Bigl(\frac{b_n^{q'}}{\mu_x}\Bigr)^{\mu_x}\prod_{x\in B_{\ssup{L}}^{\rm c}}\mu_x^{-(1-\frac {q'}q)\mu_x}.
\end{equation}
Since $q'<q$, we have that $r\equiv 1-\frac {q'}q$ is positive. According to the definition of $B_{\ssup{L}}$, the last term in \eqref{betaexpression} can be estimated by
\begin{equation}\label{lasttermesti}
\prod_{x\in B_{\ssup{L}}^{\rm c}}\mu_x^{-(1-\frac {q'}q)\mu_x}\leq \prod_{x\in B_{\ssup{L}}^{\rm c}}\bigl(L^{-r}b_n^{(q'-q)}\bigr)^{\mu_x}.
\end{equation}

The penultimate term  in \eqref{betaexpression} can be estimated as
\begin{equation}\label{second-last}
\prod_{x\in B_{\ssup{L}}}\Bigl(\frac{b_n^{q'}}{\mu_x}\Bigr)^{\mu_x}\leq C^k\prod_{x\in B_{\ssup{L}}}b_n^{(q'-q)\mu_x},
\end{equation}
since we have, using also Jensen's inequality for the logarithm,
$$
\begin{aligned}
\prod_{x\in B_{\ssup{L}}}\Bigl(\frac{b_n^{q'}}{\mu_x}\Bigr)^{\mu_x}
&=\exp\Bigl\{\Bigl(\sum_{y\in B_{\ssup{L}}}\mu_y\Bigr)\sum_{x\in B_{\ssup{L}}}\frac{\mu_x}{\sum_{y\in B_{\ssup{L}}}\mu_y}\log\frac{b_n^{q'}}{\mu_x}\Bigr\}\\
&\leq\exp\Bigl\{\Bigl(\sum_{y\in B_{\ssup{L}}}\mu_y\Bigr)\log\sum_{x\in B_{\ssup{L}}}\frac{b_n^{q'}}{\sum_{y\in B_{\ssup{L}}}\mu_y}\Bigr\}\\
&=\prod_{x\in B_{\ssup{L}}}\Bigl(\frac{b_n^{q'}\#B_{\ssup{L}}}{\sum_{y\in B_{\ssup{L}}}\mu_y}\Bigr)^{\mu_x}.
\end{aligned}
$$
Now use that $\#B_{\ssup{L}}\leq\#B\leq C\alpha_n^d=Ckb_n^{-q}$ and observe that there is a constant $C>0$ such that $(\frac kl)^l\leq C^k$, for any $l\in \{1,\dots,k\}$, since the map $y\mapsto y\log y$ is bounded on $(0,1]$. Using \eqref{lasttermesti} and \eqref{second-last} in \eqref{betaexpression}, we obtain, for some constant $C_M>0$, satisfying $\lim_{M\to\infty} C_M=0$,
\begin{equation}\label{Kesti}
\begin{aligned}
\sum_{z_1,\dots,z_k\in B}K(\mu)^{q'}&\leq C^k k! b_n^{(q'-q)k}\sum_{\mu\in A_k}\prod_{x\in B_{\ssup{L}}}C_{M,L}^{q'\mu_x}\prod_{x\in B_{\ssup{L}}^{\rm c}}L^{-r\mu_x}\\
&\leq C^k k! b_n^{(q'-q)k}\# A_k \Big(\max\{C_{M,L}^{q'}, L^{-r}\}\Big)^k\leq
C_M^{q' k}k! b_n^{(q'-q)k},
\end{aligned}
\end{equation}
where we choose $L$ in dependence on $M$ such that $\lim_{M\to\infty}\max\{C_{M,L}^{q'}, L^{-r}\}=0$, and we estimated $\# A_k =\binom {k+|B|}{|B|}\leq e^{o(k)}$ (recall that $k=n\alpha_n^2$).

Using \eqref{Kesti} and \eqref{Gesti} in \eqref{cuttingesti3}, we arrive at
\begin{equation}
\E\bigl[\langle \ell_n, Y^{\ssup{>Mb_n}}\rangle^k\bigr]\leq C_M^k k!\alpha_n^2k \alpha_n^{-\frac 1{q'}dk}\Bigl(k! b_n^{(q'-q)k}\Bigr)^{\frac 1{q'}}.
\end{equation}
Now recall that $b_n^q\alpha_n^d=k=n\alpha_n^{-2}$ and use Stirling's formula to see that the right hand side of this estimate is bounded from above by $C_M^k (nb_n)^k$ for some $C_M\to 0$ as $M\to\infty$. This ends the proof of Proposition~\ref{fieldcut}.
\end{Proof}
\qed

\end{Proof}
\qed

\subsection{Smoothing the scenery}\label{sec-smooth} 

\noindent In this section we provide the smoothing argument for the field. This
will be an adaptation of results of \cite[Sect.~3]{AC02}. Fix some smooth, rotationally invariant, and
$L^1$-normalized function $\kappa\colon\R^d\to[0,\infty)$ with $\supp(\kappa)\subset Q_1$, and put $\kappa_\delta(\cdot)
=\delta^{-d}\kappa(\cdot/\delta)$ for some small $\delta>0$. 
The convolution of two functions $f,g\colon 
\R^d\to\R$ is denoted by $f*g$. Assume any of the cases (V) and (L) and 
choose $(b_n)_n$ and $(\alpha_n)_n$ 
according to \eqref{bnalphachoice}.
We consider the rescaled and cut-down field
$\overline Y_n^{\ssup{\leq M}}\colon\R^d\to[-M,M]$; see \eqref{ydeco}. Recall the scaled and normalized
local times $L_n$ from \eqref{Lndef}. By $\Bcal_M(\R^d)$ we denote the set of all measurable functions $\R^d\to[-M,M]$.

\begin{lemma}[Scenery smoothing]\label{FieSmoo} Fix $M>0$. Then, for any $\eps>0$,
\begin{equation}\label{uniformsmoothing}
\lim_{\delta\downarrow 0}\limsup_{n\to\infty}\frac{\alpha_n^2}n\log \sup_{f\in \Bcal_M(\R^d)}\P\bigl(|
\langle L_n,[f-f*\kappa_\delta]\rangle|>\eps)=-\infty.
\end{equation}
In particular,
\begin{equation}\label{fieldsmoothing}
\lim_{\delta\downarrow 0}\limsup_{n\to\infty}\frac{\alpha_n^2}n\log \P\bigl(|
\langle L_n,[\overline Y_n^{\ssup{\leq M}}-\overline Y_n^{\ssup{\leq M}}*\kappa_\delta]\rangle|>\eps)=-\infty.
\end{equation}
\end{lemma}

\begin{Proof}{Proof} Certainly, it suffices to prove \eqref{uniformsmoothing} for $M=1$. We adapt the proof of \cite[Lemma~3.1]{AC02}, which is the same statement for $M=1$ and Brownian motion instead of random walk in Brownian scaling. We shall write $\Bcal$ instead of $\Bcal_1(\R^d)$.

Since all exponential moments of the steps are assumed finite, we have
$$
\lim_{R\to\infty}\limsup_{n\to\infty}\frac{\alpha_n^2}n\log \P\bigl(\supp(\ell_n)\not\subset B_{R_n}\bigr)=-\infty,
$$
where $R_n=R n\alpha_n^{-1}$. Hence, it suffices to show, for every $R>0$,
\begin{equation}
\lim_{\delta\downarrow 0}\limsup_{n\to\infty}\frac{\alpha_n^2}n\log \sup_{f\in \Bcal}\P\bigl(|
\langle L_n,f-f*\kappa_\delta\rangle|>\eps, \supp(\ell_n)\subset B_{R_n})=-\infty.
\end{equation}
We prove this only without absolute value signs, since the complementary inequality is proved in the same way. Fix $f\in\Bcal$. Chebyshev's inequality yields, for any $a>0$,
\begin{equation}\label{aimsmooth}
\begin{aligned}
\P\bigl(
\langle &L_n,f-f*\kappa_\delta\rangle>\eps, \supp(\ell_n)\subset B_{R_n})\\
&\leq \E\Bigl[\exp\Bigl\{a \frac n{\alpha_n^{2}}\langle L_n,f-f*\kappa_\delta\rangle\Bigr\}\1\{\supp(\ell_n)\subset B_{R_n}\}\Bigr]e^{-a\eps n\alpha_n^{-2}}.
\end{aligned}
\end{equation}
Introduce a discrete version $\varphi_n\colon \Z^d\to\R$ of $f-f*\kappa_\delta$ by
\begin{equation}
\varphi_n(z)=\alpha_n^{d}\int_{z\alpha_n^{-1}+[0,\alpha_n^{-1})^d}[f-f*\kappa_\delta](x)\,\d 
x,\qquad z\in\Z^d.
\end{equation}
Note that 
\begin{equation}\label{phinrescale}
\begin{aligned}
\frac n{\alpha_n^2} \langle L_n,f-f*\kappa_\delta\rangle&=\alpha_n^{d-2}\int[f-f*\kappa_\delta](x) 
\ell_n\bigl(\lfloor x\alpha_n\rfloor\bigr)\,\d x
=\alpha_n^{-2}\sum_{z\in \Z^d}\ell_n(z)\varphi_n(z)\\
&=\alpha_n^{-2}\sum_{k=0}^n\varphi_n(S_k).
\end{aligned}
\end{equation}

We first express the expectation on the right side of \eqref{aimsmooth} in terms of an 
expansion with respect to an appropriate orthonormal system of eigenvalues 
and eigenfunctions in $\R^{B_{R_n}}$. We write $\E_z$ for 
expectation with respect to the random walk when started at $z\in\Z^d$, in 
particular $\E=\E_0$. By \eqref{phinrescale}, for any $z,\widetilde 
z\in B_{R_n}$,
\begin{equation}\label{rewriteexp}
\E_z\Bigl[\exp\Bigl\{a\frac n{\alpha_n^2} 
\langle L_n,f-f*\kappa_\delta\rangle\Bigr\}\1\{\supp(\ell_n)\subset B_{R_n}\}\1\{S_n=\widetilde z\}\Bigr]= e^{\frac a2 \alpha_n^{-2}(\varphi_n(z)+\varphi_n(\widetilde 
z))}A^n(z,\widetilde z),
\end{equation}
where $A^n$ is the $n$-th power of the symmetric matrix $A$ 
having components
\begin{equation}
A(z,\widetilde z)=e^{\frac a2 
\alpha_n^{-2}\varphi_n(z)}p(z,\widetilde z)e^{\frac a2 \alpha_n^{-2}\varphi_n(\widetilde z)},\qquad 
z,\widetilde z\in B_{R_n}.
\end{equation}
Using an expansion in terms of the eigenvalues $\lambda_{k}(n)$, $k\in\{1,\dots,|B_{R_n}|\}$, of 
$A$ and an orthonormal basis  of $\R^{B_{R_n}}$ consisting 
of corresponding eigenfunctions $v_{k,n}$ we obtain, for any $z,\widetilde z\in B_{R_n}$, 
\begin{equation}\label{rewriteexp2}
A^n(z,\widetilde z)=\sum_{k=1}^{|B_{R_n}|}\lambda_{k}(n)^n v_{k,n}(z)v_{k,n}(\widetilde z).
\end{equation}
We assume that the 
eigenvalues $\lambda_{k}(n)$ are in decreasing order, and the principal 
eigenvector $v_{1,n}$ is positive in $B_{R_n}$.

Now we use this for the expectation on the right side of \eqref{aimsmooth}, which is equal to the 
sum over $\widetilde z\in B_{R_n}$ of the left side of 
\eqref{rewriteexp} at $z=0$. We obtain an upper bound by summing the right hand 
side of \eqref{rewriteexp2} over $z,\widetilde z\in B_{R_n}$. Continuing the 
upper bound with the help of Parseval's identity gives
\begin{equation}\label{lowupp}
\begin{aligned}
\E\Bigl[&\exp\Bigl\{a\frac n{\alpha_n^2} 
\langle L_n, f-f*\kappa_\delta\rangle\Bigr\}\1\{\supp(\ell_n)\subset B_{R_n}\}\Bigr]\\
&\leq (1+o(1))\sum_{k=1}^{|B_{R_n}|}\lambda_{k}(n)^n \sum_{z,\widetilde 
z\in B_{R_n}}v_{k,n}(z)v_{k,n}(\widetilde z)
\leq (1+o(1))\lambda_{1}(n)^n\sum_{k=1}^{|B_{R_n}|}\langle 
v_{k,n},\1\rangle^2\\
&\leq (1+o(1))\lambda_{1}(n)^n|B_{R_n}|,
\end{aligned}
\end{equation}
where we denote by $\langle\cdot,\cdot\rangle$ and $\|\cdot\|_2$ the 
inner product and Euclidean norm on $\R^{B_{R_n}}$.
Recall that $R_n=Rn\alpha_n^{-1}$. Our assumptions on $(\alpha_n)_n$ imply that 
$|B_{R_n}|=e^{o(n\alpha_n^{-2})}$ as $n\to\infty$. Hence, as $n\to\infty$,
\begin{equation}\label{Evalueesti}
\frac{\alpha_n^2}n\log\E\Bigl[\exp\Bigl\{a\frac n{\alpha_n^2} 
\langle L_n, f-f*\kappa_\delta\rangle\Bigr\}\1\{\supp(\ell_n)\subset B_{R_n}\}\Bigr]\leq o(1)+\alpha_n^2 \bigl[\lambda_{1}(n)-1\bigr].
\end{equation}

Recall the Rayleigh-Ritz principle, $\lambda_{1}(n)=\max_{\|g\|\leq 1}\langle A 
g,g\rangle$, where the maximum runs over all $\ell^2$-normalized vectors $g\colon 
\Z^d\to (0,\infty)$ with support in $B_{R_n}$. Recall that $|\varphi_n|\leq 2$. Then, as 
$n\to\infty$, we have, for any $\ell^2$-normalized vector $g$,
\begin{equation}\label{eigenvesti1}
\begin{aligned}
\alpha_n^2 \bigl[\langle A g,g\rangle-1\bigr]&=\alpha_n^2 
\Bigl(\sum_{z,\widetilde z}\bigl(e^{\frac 
a2\alpha_n^{-2}[\varphi_n(z)+\varphi_n(\widetilde z)]}-1\bigr)p(z,\widetilde z)g(z)g(\widetilde 
z)+\sum_{z,\widetilde z}\bigl(p(z,\widetilde z)-\delta_{z,\widetilde 
z}\bigr)g(z)g(\widetilde z)\Bigr)\\
&=a\langle\varphi_n,g^2\rangle+a\langle\varphi_n, g\,(pg-g)\rangle+\Ocal(\alpha_n^{-2})-\alpha_n^2 \Ical^{\ssup {\rm d}}(g^2),
\end{aligned}
\end{equation}
where we recall that the walk is assumed symmetric, and we introduced its Dirichlet form, 
\begin{equation}\label{Dirichlet}
\Ical^{\ssup {\rm d}}(g^2) =\frac 12\sum_{z,\widetilde z\in\Z^d}p(z,\widetilde z)\bigl(g(z)-g(\widetilde z)\bigr)^2,\qquad g\in\ell_2(\Z^d),
\end{equation}
and we wrote $pg(z)=\sum_{\widetilde z}p(z,\widetilde z)g(\widetilde z)$. 

The second term on the right hand side of \eqref{eigenvesti1} is estimated as follows, using that $|\varphi_n|\leq 2$.
\begin{equation}\label{eigenvesti2}
\begin{aligned}
\langle\varphi_n, g\,(pg-g)\rangle&=\frac 12\sum_{z,\widetilde z}\varphi_n(z)p(z,\widetilde z)\Big[-\bigl(g(z)-g(\widetilde z)\big)^2+\big(g(\widetilde z)-g(z)\big)\big(g(z)+g(\widetilde z)\big)\Big]\\
&\leq 2\Ical^{\ssup {\rm d}}(g^2)+\sqrt{2\Ical^{\ssup {\rm d}}(g^2)}\sqrt{\frac12\sum_{z,\widetilde z}|\varphi_n(z)|p(z,\widetilde z)(g(z)+g(\widetilde z))^2}\\
&\leq 2\Ical^{\ssup {\rm d}}(g^2)+\frac 8\eps \Ical^{\ssup {\rm d}}(g^2)+\frac\eps4,
\end{aligned}
\end{equation}
where we used the inequality $\sqrt{2ab}\leq 8a/\eps+\eps b/16$ for $a,b,\eps>0$ in the last step.

The first  term on the right hand side of \eqref{eigenvesti1} is estimated as follows. We introduce $g_n(x)=g(\lfloor x\alpha_n\rfloor)$.
\begin{equation}\label{eigenvesti3}
\begin{aligned}
\langle\varphi_n,g^2\rangle&=\alpha_n^d\int\d x\,f(x) \Big(g_n^2(x)-\int\d y\,\kappa_\delta(y) g_n^2(x+y)\Big)\\
&\leq \alpha_n^d\int\d x\, \int\d y\,\kappa_\delta(y)\big|g_n^2(x)-g_n^2(x+y)\big|\\
&\leq \alpha_n^d\int\d x\, \sqrt{\int\d y\,\kappa_\delta(y)(g_n(x)-g_n(x+y))^2}\sqrt {\int\d y\,\kappa_\delta(y)(g_n(x)+g_n(x+y))^2}\\
&\leq \frac 4{\eps}\alpha_n^d\int\d x\, \int\d y\,\kappa_\delta(y)\big(g_n(x)-g_n(x+y)\big)^2+\frac\eps8 \alpha_n^d\int\d x\,\Big(g_n^2(x)+\int\d y\,\kappa_\delta(y)g_n^2(x+y)\Big)\\
&\leq \frac 4{\eps}\alpha_n^d\int\d x\, \int\d y\,\kappa_\delta(y)\big(g_n(x)-g_n(x+y)\big)^2+\frac \eps4,
\end{aligned}
\end{equation}
where we used that $|f|\leq 1$ in the second step, H\"older's inequality in the third, and the inequality $\sqrt{2ab}\leq 4a/\eps+\eps b/8$ in the fourth step. Now pick some almost everywhere differentiable function $\psi_n\colon\R^d\to\R$ such that $\psi_n(z/\alpha_n)=\alpha_n^{d/2}g(z)$ for any $z\in\Z^d$, then a Taylor expansion gives that
\begin{equation}\label{eigenvesti4}
\begin{aligned}
\alpha_n^d\int\d x\, &\int\d y\,\kappa_\delta(y)\big(g_n(x)-g_n(x+y)\big)^2\\
&=\alpha_n^{-d}\sum_{z,\widetilde z}\Big(\psi_n({\textstyle{\frac z{\alpha_n}}})-\psi_n({\textstyle{\frac {z+\widetilde z}{\alpha_n}}})\Big)^2\int_{\widetilde z/\alpha_n+[0,1/\alpha_n]^d}\d y\,\kappa_\delta(y)\\
&=\alpha_n^{-d}\sum_{z,\widetilde z}\Big(\int_0^1\d t\, \frac{\widetilde z}{\alpha_n}\cdot\nabla\psi_n({\textstyle{\frac{z+t\widetilde z}{\alpha_n}}})\Big)^2\int_{\widetilde z/\alpha_n+[0,1/\alpha_n]^d}\d y\,\kappa_\delta(y)\\
&\leq \alpha_n^{-d} \sum_{\widetilde z}\int_{\widetilde z/\alpha_n+[0,1/\alpha_n]^d}\d y\,\kappa_\delta(y)|{\textstyle{\frac {\widetilde z}{\alpha_n}}}|^2\int_0^1\d t\,\sum_z\big|\nabla\psi_n({\textstyle{\frac{z+t\widetilde z}{\alpha_n}}})\big|^2\\
&\leq C\delta^2\|\nabla\psi_n\|_2^2\leq C\delta^2\|\Gamma^{\frac 12}\nabla\psi_n\|_2^2,
\end{aligned}
\end{equation}
where we remark that $\int\d y\,\kappa_\delta(y)|y|^2\leq C\delta^2$ for some $C>0$. 
Now we specialize the choice of $\psi_n$ to
$$
\psi_n(x)=\alpha_n^{d/2}\Big[g_n(x)+\sum_{i=1}^d\big(\alpha_n x_i-\lfloor\alpha_n x_i\rfloor\big)\Big(g\big(\lfloor\alpha_n x\rfloor+{\rm e}_i\big)-g\big(\lfloor\alpha_n x\rfloor\big)\Big)\Big],
$$
where ${\rm e}_i$ denotes the $i$-th unit vector. Then $\psi_n$ is the linear interpolation of the rescaling of $g$, and  $\partial_i\psi_n(x)=\alpha_n^{d/2+1}(g(\lfloor\alpha_n x\rfloor+{\rm e}_i)-g(\lfloor\alpha_n x\rfloor))$. Similarly to \eqref{eigenvesti4}, one derives
\begin{equation}\label{eigenvesti5}
\begin{aligned}
\alpha_n^2 \Ical^{\ssup {\rm d}}(g^2)&=\int_0^1\d t\,\int_0^1\d s\,\sum_{z\in\Z^d}p(0,z)\sum_{i,j=1}^dz_iz_j\int\d x\,
\partial_i \psi_n\Big(\frac{\lfloor\alpha_n x\rfloor+tz}{\alpha_n}\Big)
\partial_j\psi_n\Big(\frac{\lfloor\alpha_n x\rfloor+sz}{\alpha_n}\Big)\\
&=\int_0^1\d t\,\int_0^1\d s\,\sum_{z}p(0,z)\sum_{i,j=1}^dz_iz_j\int\d x\,\partial_i \psi_n(x)\partial_j\psi_n(x)\\
&=\|\Gamma^{\frac 12}\nabla\psi_n\|_2^2.
\end{aligned}
\end{equation}
Now use \eqref{eigenvesti5} in \eqref{eigenvesti4} and this in \eqref{eigenvesti3}, and substitute \eqref{eigenvesti3} and \eqref{eigenvesti2} in \eqref{eigenvesti1} to obtain, for any $a>0$, for $n$ sufficiently large and all $\ell^2$-normalized $g\in\ell^2(\Z^d)$ with support in $B_{R_n}$,
$$
\alpha_n^2 \bigl[\langle A g,g\rangle-1\bigr]\leq \frac 12a\eps -\alpha_n^2\Ical^{\ssup {\rm d}}(g^2)\Big(1-C\frac {\delta^2 a}\eps\Big),
$$
for some $C>0$ which does not depend on $n$, $g$, $\eps$ or on $a$.
Now we choose $a=\eps/(2C\delta^2)$ and obtain $\alpha_n^2 \bigl[\langle A g,g\rangle-1\bigr]\leq \frac 12a\eps$. Taking the supremum over all $g$'s considered, we obtain that $\alpha_n^2[\lambda_{1}(n)-1]\leq \frac 12a\eps$. Using this in \eqref{Evalueesti} and \eqref{Evalueesti} in \eqref{aimsmooth}, we obtain that
$$
\mbox{l.h.s.~of \eqref{fieldsmoothing}}\leq \limsup_{\delta\downarrow 0}-\frac 12 a\eps=-\lim_{\delta\downarrow 0}\frac {\eps^2}{4C\delta^2}=-\infty,
$$
and the proof is finished.
\end{Proof}\qed

\subsection{Various approximations}\label{sec-appr}

\noindent In the proofs of Theorems~\ref{inter} and \ref{lin} we shall need a couple of 
approximations to the variational formulas in \eqref{K3def} and \eqref{kdq}.
In particular, we need to show that they may be approximated by finite-space approximations
and by smoothed versions of the functions involved in the variational formula.

As in Section~\ref{sec-smooth}, by $\kappa=\kappa_1\colon\R^d\to[0,\infty)$ we denote a smooth, rotationally invariant $L^1$-normalized function, and we put $\kappa_\delta(x)=\delta^{-d}\kappa_1(x\delta^{-1})$ for $\delta>0$. Hence, $\kappa_\delta$ is a smooth approximation of the Dirac measure  at zero.

\begin{lemma}[Approximations of $K_{H}$]\label{KHappr} 
For any $u>0$, 
\begin{equation}
\limsup_{\delta\downarrow 0}\limsup_{R\to\infty}K_{H}^{\ssup{0}}(u;\delta, R)\leq K_H(u)\leq\liminf_{\delta\downarrow 0}\liminf_{R\to\infty}K_{H}^{\ssup{\rm per}}(u;\delta, R),
\end{equation}
where
\begin{eqnarray}
K_{H}^{\ssup{0}}(u;\delta, R)&=&\inf\Bigl\{
    \frac 12\|\Gamma^{\frac 12}\nabla \psi\|_2^2+\Phi_H(\psi^2*\kappa_\delta,u;R)\colon \psi\in 
    H^1(\R^d),\supp(\psi)\subset Q_R,\label{KHR0def}\\
&&\qquad\qquad\qquad\qquad\qquad\qquad \|\psi\|_2=1\Bigr\},\nonumber \\
K_{H}^{\ssup{\rm per}}(u;\delta, R)&=&\inf\Bigl\{
    \frac 12\|\Gamma^{\frac 12}\nabla_{R} \psi\|_2^2+\Phi_H(\psi^2*\kappa_\delta,u;R)\colon \psi\in 
    H^1(Q_R),\|\psi\|_2=1\Bigr\},\label{KHRperdef}
\end{eqnarray}
and $\Phi_H(\psi^2,u;R)$ is defined in \eqref{PhiRdef}.
In \eqref{KHRperdef}, $\nabla_R$ denotes the gradient on the torus $Q_R$, i.e., with periodic boundary condition.
\end{lemma}

\begin{Proof}{Proof}
Fix $\delta>0$. In the first step, we carry out the limit as $R\to\infty$ on both sides to obtain
\begin{equation}\label{Rlimtwithdelta}
\limsup_{R\to\infty}K_{H}^{\ssup{0}}(u;\delta, R)\leq K_H(u;\delta)\leq\liminf_{R\to\infty}K_{H}^{\ssup{\rm per}}(u;\delta, R),
\end{equation}
where $K_H(u;\delta)$ is defined as $K_H(u)$ in \eqref{K3def} with $\Phi_H(\psi^2,u)$ replaced by $\Phi_H(\psi^2*\kappa_\delta,u)$. The proof of \eqref{Rlimtwithdelta} follows standard patterns (see the proof of \cite[Lemma~3.7]{AC02}, e.g.) and we do not carry this out here. Hence, the only thing left to do is to show that $\lim_{\delta\downarrow 0}K_H(u;\delta)=K_H(u)$. 

Using the convexity of $H$, it is easy to derive with the help of Jensen's inequality that, for any $\gamma>0$ and any $\psi$,
$$
\int H\bigl(\gamma \psi^2*\kappa_\delta(y)\bigr)\,\d y\leq \int H\bigl(\gamma\psi^2(y)\bigr)\,\d y.
$$
As a consequence, we have $\Phi_H(\psi^2*\kappa_\delta,u)\geq \Phi_H(\psi^2,u)$ and therefore $K_{H}(u;\delta)\geq K_{H}(u)$ for any $\delta>0$.  

We argue now that $\limsup_{\delta\downarrow 0}K_{H}(u;\delta)\leq K_{H}(u)$. Indeed, fix some small $\eps>0$ and pick some bounded approximative $\eps$-minimizer for $K_{H}(u)$, i.e., a bounded function $\overline \psi\in H^1(\R^d)$ satisfying $\|\overline\psi\|_2=1$ and 
$$
\frac 12\|\Gamma^{\frac 12}\nabla \overline\psi\|_2^2+\Phi_H(\overline\psi^2,u)\leq K_{H}(u) +\eps.
$$
Using the mean-value theorem and the fact that $\|\overline\psi^2*\kappa_\delta-\overline\psi^2\|_1\to0$ as $\delta\downarrow 0$ (see \cite[Th.~2.16]{LL97}), it is elementary to show that we have $\int H\bigl(\gamma \overline\psi^2*\kappa_\delta(y)\bigr)\,\d y \to \int H\bigl(\gamma \overline\psi^2(y)\bigr)\,\d y$ as $\delta\downarrow 0$, uniformly in $\gamma$ on any compact subset of $[0,\infty)$. As a consequence, we have $\lim_{\delta\downarrow 0}\Phi_H(\overline \psi^2*\kappa_\delta,u)=\Phi_H(\overline \psi^2,u)$ and therefore
\begin{equation}\label{Kdeltaconv}
\begin{aligned}
\limsup_{\delta\downarrow 0}K_{H}(u;\delta)&\leq \frac 12\|\Gamma^{\frac 12}\nabla \overline\psi\|_2^2+\limsup_{\delta\downarrow 0}\Phi_H(\overline \psi^2*\kappa_\delta,u)=\frac 12\|\Gamma^{\frac 12}\nabla \overline\psi\|_2^2+\Phi_H(\overline\psi^2,u)\\
&\leq K_{H}(u) +\eps.
\end{aligned}
\end{equation}
Now let $\eps\downarrow 0$.
\end{Proof}\qed

Lemma~\ref{KHappr} implies the corresponding statement for the case (V):

\begin{cor}[Approximations of $K_{D,q}$]\label{KDqappr} 
Fix $D>0$ and $q>1$ and recall that $\frac 1p+\frac 1q=1$. Then
\begin{equation}
\limsup_{R\to\infty}K_{D,q}^{\ssup{0}}(R)\leq K_{D,q}\leq\liminf_{\delta\downarrow 0}\liminf_{R\to\infty}K_{D,q}^{\ssup{\rm per}}(\delta,R),
\end{equation}
where
\begin{eqnarray}
K_{D,q}^{\ssup{0}}(R)&=&\inf\Bigl\{
    \frac 12\|\Gamma^{\frac 12}\nabla \psi\|_2^2+D\|\psi^2\|_{p}^{-q}\colon\psi\in H^1(\R^d),\supp(\psi)\subset Q_R,\|\psi\|_2=1\Bigr\},\\
K_{D,q}^{\ssup{\rm per}}(\delta,R)&=&\inf\Bigl\{
    \frac 12\|\Gamma^{\frac 12}\nabla_R \psi\|_2^2+D\|\psi^2*\kappa_\delta\|_{p}^{-q}\colon\psi\in H^1(Q_R),\|\psi\|_2=1\Bigr\},
\end{eqnarray}
and $\nabla_R$ is the  gradient on the torus $Q_R$, i.e., with periodic boundary condition.
\end{cor}

\begin{Proof}{Proof}
We apply Lemma~\ref{KHappr} to the special choice $u=1$ and $H(t)=\widetilde D t^p$, where $p$ and $\widetilde D$ are as in \eqref{cumgenfct}. It is easy to see that for this choice of $H$, we have $\Phi_H(\psi^2,1)=D\|\psi^2\|_p^{-q}$. 
\end{Proof}\qed

\section{Proof of the upper bounds in Theorems~\ref{inter} and \ref{lin}}
\label{s:ub}

\noindent This section is devoted to the proof of the upper bounds in
Theorems~\ref{inter} and \ref{lin}. They are in Sections~\ref{subs:inter} and \ref{subs:ld}, respectively. Our proofs essentially follow the proof of \cite[Theorem~2.2]{AC02}.

\subsection{Very-large deviation case (Theorem~\ref{inter})}
\label{subs:inter}

\noindent In this section we are  under Assumption (Y) with $q>\frac d2$, and consider a sequence $(b_n)$ with $1\ll b_n \ll
n^{\frac 1q}$. We have to smoothen the scenery, as we have explained at the beginning of Section~\ref{sec-prep}. In order to do this, we have to cut down the scenery to bounded size. As soon as the smoothing argument has been carried out, we may relax the boundedness assumption.

Recall the scaled and normalized local times $L_n$ from \eqref{Lndef}
and the scaled normalized scenery $\overline Y_n$ from \eqref{Yscaled}.
Recall the notation $y^{\ssup{\leq M}}=[y\wedge M]\vee(-M)$ from \eqref{ydeco}, and recall the
delta-approximation $\kappa_\delta\colon\R^d\to[0,\infty)$ to the Dirac measure from the beginning
of Section~\ref{sec-smooth}.

Note that, for any $M,\eps,\delta>0$,
\begin{equation}
\begin{aligned}
\P(\smfrac 1n Z_n>b_n)&\leq \P\bigl(\langle \overline Y_n^{\ssup{\leq M}}*\kappa_\delta,
L_n\rangle >1-2\eps\bigr)\\
&\qquad+\P\bigl(\langle |\overline Y_n^{\ssup{\leq M}}*\kappa_\delta
-\overline Y_n^{\ssup{\leq M}}|,L_n\rangle>\eps\bigr)
+\P(\langle \overline Y_n^{\ssup{>M}},L_n\rangle>\eps).
\end{aligned}
\end{equation}
Recall that, by our choice of $\alpha_n$, we have
\begin{equation}\label{alphachoiceI}
n^{\frac d{d+2}}b_n^{\frac {2q}{d+2}}=\frac n{\alpha_n^2}.
\end{equation}
Hence, by Proposition~\ref{fieldcut}, Lemma~\ref{FieSmoo} and Corollary~\ref{KDqappr}, it suffices to
prove, for any $M,\delta>0$ and $R\in\N$,
\begin{equation}\label{upperaim}
\limsup_{\eps\downarrow 0}\limsup_{n\to\infty}\frac {\alpha_n^2}n\log 
\P\bigl(\langle \overline Y_n^{\ssup{\leq M}}*\kappa_\delta,L_n\rangle >1-2\eps\bigr)\leq -K_{D,q}^{\ssup{\rm per}}(\delta,R),
\end{equation}
where $K_{D,q}^{\ssup{\rm per}}(\delta,R)$ is defined in Corollary~\ref{KDqappr}. Note that
$$
\begin{aligned}
\langle \overline Y_n^{\ssup{\leq M}}*\kappa_\delta,L_n\rangle
 &=\frac 1{b_n\alpha_n^d}\sum_{z\in\Z^d}\bigl[\bigl(Y(z)\wedge(Mb_n)\bigr)\vee 
(-Mb_n)\bigr]L_n*\kappa_\delta\Bigl(\frac z{\alpha_n}\Bigr)\\
&\leq \frac 1{b_n\alpha_n^d}\sum_{z\in\Z^d}\bigl[Y(z)\vee 
(-Mb_n)\bigr]L_n*\kappa_\delta\Bigl(\frac z{\alpha_n}\Bigr).
\end{aligned}
$$
Introduce the cumulant generating function of $Y(0)\vee (-M)$, 
$$
H_M(t)=\log \E[e^{t[Y(0)\vee (-M)]}].
$$
Using the exponential Chebyshev inequality and carrying out the expectation over the scenery, we obtain, for any $\gamma>0$,  the upper bound
\begin{equation}\label{Iproof1}
\begin{aligned}
\P\bigl(\langle \overline Y_n^{\ssup{\leq M}}*\kappa_\delta,L_n\rangle >1-2\eps\bigr)
&\leq\E\Bigl[ e^{-\gamma(1-2\eps) nb_n} e^{\gamma n\alpha_n^{-d}\sum_z[Y(z)\vee (-Mb_n)]L_n*\kappa_\delta\bigl(\frac z{\alpha_n}\bigr)}\Bigr]\\
&\leq \E\Bigl[e^{-\gamma (1-2\eps) nb_n}\exp\Bigl\{\sum_{z\in\Z^d} H_{Mb_n}
\Bigl(\gamma n\alpha_n^{-d}L_n*\kappa_\delta\Bigl(\frac z{\alpha_n}\Bigr)\Bigr)\Bigr\}\Bigr].
\end{aligned}
\end{equation}
Since $H_{Mb_n}$ is convex and satisfies $H_{Mb_n}(0)=0$, it is also 
superadditive. 
Hence,
for any $\gamma>0$ and any $x\in\Z^d$, we have
\begin{equation}\label{Hperesti}
\sum_{k\in \Z^d} H_{Mb_n}\bigl(\gamma \ell_n(x+2k\lfloor R\rfloor)\bigr)\leq 
H_{Mb_n}(\gamma\ell_n^{\ssup{R}}(x)),
\end{equation}
Therefore, the right hand in (\ref{Iproof1}) side does not become smaller if $L_n*\kappa_\delta$ is replaced by its periodized version, $(L_n*\kappa_\delta)^{\ssup{R}}(x)=\sum_{k\in\Z^d}L_n*\kappa_\delta(x+kR)$, for $x\in [-R,R]^d$. Furthermore, note that 
$$
(L_n*\kappa_\delta)^{\ssup{R}}(x)=\sum_{k\in\Z^d}\int_{\R^d} L_n(y)\kappa_\delta(x+kR-y)\,\d y=\int_{\R^d} L_n^{\ssup{R}}(y)\kappa_\delta(x-y)\,\d y=L_n^{\ssup{R}}*\kappa_\delta(x),
$$
for any $x\in[-R,R]^d$. Hence, we may replace $L_n$ on the right of \eqref{Iproof1} by its periodized version $L_n^{\ssup{R}}$. 

According to \eqref{cumgenfct}, for any $\eps>0$, we may choose a $c(\eps)>0$
such that 
\begin{equation}\label{Hesti}
H(t)\leq c(\eps) t+(1+\eps)\widetilde D\, t^p,\qquad t\in[0,\infty).
\end{equation}
Since $e^{H_M(t)}\leq e^{H(t)}+1$, we also have the estimate in \eqref{Hesti} for $H_{Mb_n}$ instead of $H$.  Hence, since $\kappa_\delta$ and $L_n$ are $L^1$-normalized,
\begin{equation}\label{Intercalc1}
\P\bigl(\langle \overline Y_n^{\ssup{\leq M}}*\kappa_\delta,L_n\rangle >1-2\eps\bigr)
\leq\E\Bigl[e^{-\gamma(1-2\eps) nb_n}e^{c(\eps)\gamma n}\exp\Bigl\{\gamma^p(\widetilde D+\eps)\alpha_n^d (n\alpha_n^{-d})^p\| L_n^{\ssup{R}}*\kappa_\delta\|_p^p\Bigr\}\Bigr].
\end{equation} 
We choose the value of $\gamma$ optimal for $\eps=0$, which is 
\begin{equation}\label{gammachoiceI}
\gamma=\frac{\alpha_n^d}n b_n^{\frac 1{p-1}}\Bigl(p\widetilde D\|L_n*\kappa_\delta^p\|_p^p\Bigr)^{-\frac 1{p-1}}=\frac{\alpha_n^d}n b_n^{\frac 1{p-1}}Dq\bigl\|L_n^{\ssup{R}}*\kappa_\delta\bigr\|_p^{-q},
\end{equation}
where we recalled that $1=\frac 1p+\frac 1q$ and $\widetilde D=(q-1)(Dq^q)^{\frac 1{1-q}}$. Note that the map $\mu\mapsto \|\mu*\kappa_\delta\|_p$ is bounded and continuous (in the weak $L^1$-topology)  on the set of probability measures on $[-R,R]^d$. Indeed, the continuity is seen with the help of Lebesgue's theorem, and the boundedness follows from the following application of Jensen's inequality:
\begin{equation}\label{Jensentrick}
\begin{aligned}
\|\mu*\kappa_\delta\|_p^p&=(2R)^d\int_{[-R,R]^d} \frac{\d x}{(2R)^d}\Big|\int_{\R^d}\mu(\d y)\kappa_\delta(x-y)\Big|^p\\
&\geq (2R)^d \Bigl(\int_{[-R,R]^d}\frac{\d x}{(2R)^d}\int_{\R^d}\mu(\d y)\kappa_\delta(x-y)\Bigr)^p\\
&=(2R)^{d(1-p)},
\end{aligned}
\end{equation}
since $\kappa_\delta$ is $L^1$-normalized.

Recall that $b_n^q=n\alpha_n^{-(d+2)}$. For the choice of $\gamma$ in \eqref{gammachoiceI}, for large $n$, we can estimate the first two terms in the expectation on the right of \eqref{Intercalc1} by $e^{-\gamma(1-2\eps) nb_n}e^{c(\eps)\gamma n}\leq e^{-\gamma(1-3\eps) nb_n}$, since we have in particular  $\gamma\ll b_n$.

Substituting $\gamma$ in \eqref{Intercalc1}, we obtain
\begin{equation}\label{Interdev1}
\P\bigl(\langle \overline Y_n^{\ssup{\leq M}}*\kappa_\delta,L_n\rangle >1-2\eps\bigr)
\leq\E\Bigl[\exp\Bigl\{-(D+\eps C)\frac n{\alpha_n^2}\bigl\|L_n^{\ssup{R}}*\kappa_\delta\bigr\|_p^{-q}\Bigr\}\Bigr],
\end{equation}
where $C>0$ depends on $D, R$ and $q$ only. Now we can finally apply the large deviation principle in Lemma~\ref{LDP}(ii) to the right hand side of \eqref{Interdev1}. This yields the estimate in \eqref{upperaim} without $\limsup_{\eps\downarrow 0}$ and with $D$ replaced by $D+\eps C$. Letting $\eps\downarrow 0$, we easily see that \eqref{upperaim} is satisfied, which ends the proof of the upper bound in Theorem~\ref{inter}.

\subsection{Large-deviation case (Theorem~\ref{lin})}
\label{subs:ld}

\noindent In this section, we prove the upper bound in Theorem~\ref{lin}, i.e., in the case (L). The proof follows the pattern of the corresponding proof in \cite{AC02} and is analogous to the proof of Theorem~\ref{inter} in Section~\ref{subs:inter}, and hence we keep it short. Pick $b_n=1$ and $\alpha_n=n^{\frac 1{d+2}}$,
in accordance with \eqref{bnalphachoice}. Furthermore, fix $u>0$.

By Proposition~\ref{fieldcut} and Lemmas~\ref{FieSmoo} and \ref{KHappr}, it is sufficient to prove that, for any $\delta>0$ and $R\in\N$,
\begin{equation}\label{upperaimlar}
\limsup_{\eps\downarrow0}\limsup_{n\to\infty}n^{-\frac d{d+2}}\log 
\P\bigl(\langle \overline Y_n^{\ssup{\leq M}}*\kappa_\delta,L_n\rangle >u-\eps\bigr)\leq -K_{H}^{\ssup{\rm per}}(u;\delta, R),
\end{equation}
where $K_{H}^{\ssup{\rm per}}(u;\delta,R)$ is defined in Lemma~\ref{KHappr}. Fix a small $\eps>0$. Analogously to  \eqref{Iproof1}, we have the estimate
\begin{equation}
\P\bigl(\langle \overline Y_n^{\ssup{\leq M}}*\kappa_\delta,L_n\rangle >u-\eps \bigr)\leq \E\Bigl[e^{-\gamma (u-2\eps) n}\exp\Bigl\{\sum_{z\in\Z^d} H_{M}
\Bigl(\gamma n\alpha_n^{-d}L_n^{\ssup{R}}*\kappa_\delta\Bigl(\frac z{\alpha_n}\Bigr)\Bigr)\Bigr\}\Bigr],
\end{equation}
for any $\gamma>0$. Replacing $\gamma n\alpha_n^{-d}$ by $\gamma$, turning the sum into an integral, passing to the optimum over $\gamma$ and using the notation in \eqref{PhiRdef}, we obtain
\begin{equation}\label{Lcalc}
\P\bigl(\langle \overline Y_n^{\ssup{\leq M}}*\kappa_\delta,L_n\rangle >u\bigr)\leq\E\Bigl[\exp\Bigl\{-\frac n{\alpha_n^2}\Phi_{H_M}(L_n^{\ssup{R}}*\kappa_\delta,u-2\eps;R)\Bigr\}\Bigr],
\end{equation}
where we also recall that $\alpha_n^d=n\alpha_n^{-2}$. Again, for fixed $\delta>0$ and $R>0$, we can let $M\to\infty$ and $\eps\downarrow 0$ to replace $\Phi_{H_M}(L_n^{\ssup{R}}*\kappa_\delta,u-2\eps;R)$ by $\Phi_{H}(L_n^{\ssup{R}}*\kappa_\delta,u;R)$ on the right side of \eqref{Lcalc}. Analogously to \eqref{Jensentrick}, one shows that $\Phi_H(\psi^2,u)\leq |Q_R| \, \sup_{\gamma>0}\bigl(\gamma u -H(\gamma)\bigr)<\infty$ for any continuous $\psi\colon Q_R\to [0,\infty)$ satisfying $\int_{Q_R}\psi^2=1$. Hence, the map $\mu\mapsto \Phi_{H}(\mu*\kappa_\delta,u;R)$ is bounded and continuous on the set of probability measures on $Q_R$, and we may apply the large deviation principle in Lemma~\ref{LDP}(ii). This, followed by $\eps\downarrow 0$, implies that \eqref{upperaimlar} holds for any $\delta>0$ and $R\in\N$. This finishes the proof of the upper bound in Theorem~\ref{lin}.

\section{Proofs of the lower bounds in Theorems~\ref{inter} and \ref{lin}}
\label{s:lb}

In this section we prove the lower bounds in Theorems~\ref{inter} and \ref{lin}. Our proofs are variants of the analogous proofs in \cite{AC02}; they roughly follow the heuristics in Section~\ref{sec-heur}.

\subsection{Very-large deviation case (Theorem~\ref{inter})}\label{sec-lowinter}

\noindent Suppose we are in the case (V) and pick sequences $(b_n)_n$ and $(\alpha_n)_n$ 
as in \eqref{bnalphachoice}. Fix $R>0$ and any continuous positive function $\varphi\colon Q_R\to(0,\infty)$. 
Recall the scaled local times and scenery, $L_n$ and $\overline Y_n$, in \eqref{Lndef} and 
\eqref{Yscaled}. 

If $\overline Y_n\geq 
\varphi$ on $
Q_R$ and $\supp(L_n)\subset Q_R$, then
\begin{equation}
Z_n=b_n n \langle L_n, \overline Y_n\rangle
\geq b_n n  \langle L_n, \varphi\rangle.
\end{equation}
Hence, we obtain the lower bound, for any $n\in\N$,
\begin{equation}\label{interlow1}
\P\bigl(\smfrac 1n Z_n>b_n\bigr)\geq \P(\langle L_n,\varphi\rangle \geq 
1, \supp(L_n)\subset Q_R)\,\P\bigl(\overline Y_n\geq 
\varphi\mbox{ on }Q_R\bigr).
\end{equation}

With the help of the large deviation principle in Lemma~\ref{LDP}(i) it is easy to deduce that
\begin{equation}\label{interlow2}
\begin{aligned}
\lim_{n\to\infty}&\frac{\alpha_n^2}n\log \P\bigl(\langle L_n,\varphi\rangle 
\geq 1, \supp(L_n)\subset Q_R\bigr)\\
&=-\inf\bigl\{\Ical_R(\psi^2)\colon \psi\in H^1(\R^d),\supp(\psi)\subset Q_R,\|\psi\|_2=1, \langle 
\psi^2,\varphi\rangle\geq 1 \bigr\}.
\end{aligned}
\end{equation}
{From} Lemma~\ref{LDPfieldcont} we have, recalling that $n\alpha_n^{-2}=\alpha_n^d b_n^q$,
\begin{equation}\label{Ylowinter}
\liminf_{n\to\infty}\frac{\alpha_n^2}n\log \P\bigl(\overline Y_n\geq 
\varphi\mbox{ on }Q_R\bigr)\geq -D\|\varphi\|_q^q.
\end{equation}
Using \eqref{interlow2} and \eqref{Ylowinter} in \eqref{interlow1} and optimizing on $\varphi$, we obtain the lower bound
\begin{equation}\label{lowPinter}
\liminf_{n\to\infty}\frac{\alpha_n^2}n\log 
\P\bigl(\smfrac 1n Z_n>b_n\bigr)\geq -\widetilde K_{D,q}^{\ssup{0}}(R),
\end{equation}
where
\begin{equation}
\widetilde K_{D,q}(R)=\inf_{\psi\in H^1(\R^d)\colon \|\psi\|_2=1,
\supp(\psi)\subset B_R}\Bigl(\Ical_R(\psi^2)+D 
\inf_{\varphi\in\Ccal_+(Q_R)\colon \langle \psi^2,\varphi\rangle \geq 1}
\|\varphi\|_q^q\Bigr).
\end{equation}
It is easy to see that the inner infimum is equal to $\|\psi^2\|_p^{-q}$. Hence, $\widetilde K_{D,p}(R)=K_{D,p}^{\ssup{0}}(R)$ as defined in Corollary~\ref{KDqappr}. Now  Corollary~\ref{KDqappr} finishes the proof 
of the lower bound in Theorem~\ref{inter}.

\subsection{Large-deviation case (Theorem~\ref{lin})}

\noindent Recall from Section~\ref{sec-heur} that $\frac 1n Z_n=\langle L_n, \overline Y_n\rangle$. We want to apply the large deviation principles of Lemma~\ref{LDP}(i) for $L_n$ and Lemma~\ref{LDPfieldL} for $\overline Y_n$. However, as has been pointed out in \cite{AC02}, the map $(\mu,f)\mapsto \langle \mu,f\rangle$ is not continuous in the product of the weak topologies. Hence, we partially follow the strategy of \cite{AC02} and use Lemma~\ref{FieSmoo} to smoothen the field $\overline Y_n$. In order to apply Lemma~\ref{FieSmoo}, we first have to cut down the field to bounded size, which we do with the help of Proposition~\ref{fieldcut}. However, this works only for cutting the {\it large\/} values of the field, but not the small ones. In order to be able to use also a lower bound for the field, we intersect with the event that $Y(z)\geq -M$ for all $z$'s appearing, and use a large deviation principle for the conditional field. 

Let us turn to the details. Let $u>0$ satisfying $u\in\supp(Y(0))^\circ$. We fix small parameter $\eps,\delta>0$ such that $u+\eps\in\supp(Y(0))^\circ$ and large parameters $M$ and $R$. On the intersection of the events $\{\supp(L_n)\subset Q_R\}$ and $\{Y(z)\geq -M\, \forall z\in B_{R\alpha_n}\}$, we can estimate
$$
\frac 1n Z_n=\langle L_n, \overline Y_n\rangle \geq\langle L_n, \overline Y_n^{\ssup{\leq M}}\rangle=\langle L_n*\kappa_\delta, \overline Y_n^{\ssup{\leq M}}\rangle+\langle L_n, \overline Y_n^{\ssup{\leq M}}-\overline Y_n^{\ssup{\leq M}}*\kappa_\delta\rangle.
$$
We write $\P^{\ssup{>-M}}$ for the conditional measure $\P(\,\cdot\,|\,Y(z)\geq -M\, \forall z\in\Z^d)$. Hence, we obtain the lower bound
\begin{equation}
\begin{aligned}
\P(\smfrac 1n Z_n>u)&\geq \P^{\ssup{>-M}}\Bigl(\supp(L_n)\subset Q_R, \langle L_n*\kappa_\delta, \overline Y_n^{\ssup{\leq M}}\rangle>u+\eps\Bigr)\P(Y(0)\geq -M)^{|B_{R\alpha_n}|}\\
&\qquad-\P(\langle L_n, \overline Y_n^{\ssup{\leq M}}-\overline Y_n^{\ssup{\leq M}}*\kappa_\delta\rangle>\eps).
\end{aligned}
\end{equation}
Using  Lemma~\ref{FieSmoo} for the last term on the right hand side, and noting that $\P(Y(0)\geq -M)\to0$ as $M\to\infty$, it becomes clear that it suffices to estimate the first term on the right side. In order to do this, fix a positive continuous function $g\colon Q_R\to(0,\infty)$ satisfying $\int_{Q_R}g(x)\,\d x =1$ such that $g$ can be extended to an element of $H^1(\R^d)$. Let $B_\eps(g)$ denote a weak $\eps$-neighborhood of $g$. Then we have
$$
\begin{aligned}
\P^{\ssup{>-M}}\Bigl(&\supp(L_n)\subset Q_R, \langle L_n*\kappa_\delta, \overline Y_n^{\ssup{\leq M}}\rangle>u+\eps\Bigr)\\
&\geq \P(L_n\in B_\eps(g), \supp(L_n)\subset Q_R)\P^{\ssup{>-M}}\Bigl(\langle g*\kappa_\delta, \overline Y_n^{\ssup{\leq M}}\rangle>u+2\eps\Bigr).
\end{aligned}
$$
According to Lemma~\ref{LDP}, the first term on the right is equal to $\exp\{-n\alpha_n^{-2} \inf_{\psi^2\in B_\eps(g)}\Ical_R(\psi^2)(1+o(1))\}$, and according to Lemma~\ref{LDPfieldL}, the latter term is equal to $\exp\{-n\alpha_n^{-2} \Phi_{\widetilde H_M}(g*\kappa_\delta,u-2\eps,R)(1+o(1))\}$. Summarizing, we obtain, for any $R>0$ and any continuous positive function $g\colon Q_R\to(0,\infty)$, if $M$ is sufficiently large and $\delta>0$ sufficiently small,
\begin{equation}
\liminf_{n\to\infty}\frac {\alpha_n^2}n\log \P(\smfrac 1n Z_n>u)\geq -\Bigl[\Ical_R(g)+\Phi_{\widetilde H_M}(g*\kappa_\delta,u+2\eps,R)\Bigr]+\eta_M,
\end{equation}
for some $\eta_M\downarrow 0$ as $M\to\infty$. Passing to the infimum over all $g$ and writing $\psi^2$ instead of $g$, we obtain
\begin{equation}
\liminf_{n\to\infty}\frac {\alpha_n^2}n\log \P(\smfrac 1n Z_n>u)\geq -\inf\limits_{\psi\in H^1(\R^d)\colon \supp(\psi)\subset Q_R}\Bigl[\Ical_R(\psi^2)+\Phi_{\widetilde H_M}(\psi^2*\kappa_\delta,u+2\eps,R)\Bigr]+\eta_M.
\end{equation}
Since $\psi^2*\kappa_\delta$ is bounded uniformly in $\psi$, and since $\widetilde H_M(t)\to H(t)$ as $M\to\infty$, uniformly in $t$ on compacts, we can let $M\to\infty$. Furthermore, we also let $\eps\downarrow0$ and obtain
\begin{equation}
\liminf_{n\to\infty}\frac {\alpha_n^2}n\log \P(\smfrac 1n Z_n>u)\geq -K_H(u;\delta,R),
\end{equation}
for any $\delta>0$ and $R>0$, where $K_H^{\ssup{0}}(u;\delta,R)$ is defined in Lemma~\ref{KHappr}. Now use Lemma~\ref{KHappr} to finish the proof of the lower bound in Theorem~\ref{lin}.

\section{Appendix: Proof of the large deviation principle for the local times}\label{sec-proofLDP}

\noindent In this section, we prove the scaled large deviation principles in 
Lemma~\ref{LDP}.
Although the statement should be familiar to experts and the proof is 
fairly standard, we could not find it
in the literature.  Therefore, we provide a proof. Let us mention that the lower bound of the following
Lemma~\ref{lem-cumul} (without the indicator on $\{\supp(L_n)\subset Q_R\}$, however) is contained in \cite{CL02}.

Fix $R>0$.  For bounded and continuous functions $ f \colon Q_R\to\R$, we 
denote by
\begin{equation}\label{Lapeigenv}
\lambda_R( f )=\max\Bigl\{\langle f ,\psi ^2\rangle-\frac{1}{2}\|\Gamma^{\frac 
12}\nabla \psi \|_2^2\colon \psi \in H^1(\R^d),\supp( \psi )\subset Q_R,\| 
\psi \|_2  = 1\Bigr\}
\end{equation}
the principal eigenvalue of the operator 
$\frac{1}{2}\nabla\cdot\Gamma\nabla+ f $ in $Q_R$ with Dirichlet boundary
condition. (We denote the inner product and norm on $L^2(Q_R)$ by $\langle\cdot,\cdot\rangle$ 
and $\|\cdot\|_2$.)
The main step in the proof of Lemma~\ref{LDP}(i) is the 
following.

\begin{lemma}\label{lem-cumul}
For any bounded and continuous function $ f \colon Q_R\to\R$, the 
limit
\begin{equation}\label{cumul}
\lim_{n\to\infty}\frac {\alpha_n^2}n\log \E\Bigl[\exp\Bigl\{\frac 
n{\alpha_n^2} \langle f ,L_n\rangle\Bigr\}\1\{\supp(L_n)\subset Q_R\}\Bigr]
\end{equation}
exists and is equal to $\lambda_R( f )$.
\end{lemma}

\begin{Proof}{Proof} In the following, we abreviate $B=B_{R\alpha_n}$.
Introduce a scaled version $ f _n\colon \Z^d\to\R$ of $ f $ by
\begin{equation}
 f _n(z)=\alpha_n^{d}\int_{z\alpha_n^{-1}+[0,\alpha_n^{-1})^d} f (x)\,\d 
x,\qquad z\in\Z^d.
\end{equation}
Note that $ f _n(\lfloor \cdot\,\alpha_n\rfloor)\to  f $ 
uniformly on $Q_R$.  Furthermore, note that 
\begin{equation}
\frac n{\alpha_n^2} \langle f ,L_n\rangle=\alpha_n^{d-2}\int_{Q_R} f (x) 
\ell_n\bigl(\lfloor x\alpha_n\rfloor\bigr)\,\d x
=\alpha_n^{-2}\sum_{z\in B}\ell_n(z) f _n(z)
=\sum_{k=0}^{n-1}\alpha_n^{-2} f _n(S_k).
\end{equation}

For notational convenience, we assume that $\alpha_n^2$ and $n\alpha_n^{-2}$ are integers. Using the Markov property, we split the expectation over the path $(S_0,\dots,S_n)$ into $n\alpha_n^{-2}$ expectations over paths of length $\alpha_n^2$. By $\E_z$ we denote the expectation with respect to the random walk starting at $z\in\Z^d$, then we have
\begin{equation}\label{expansionexp}
\begin{aligned}
\E\Bigl[&\exp\Bigl\{\frac 
n{\alpha_n^2} \langle f ,L_n\rangle\Bigr\}\1\{\supp(L_n)\subset Q_R\}\Bigr]\\
&=\E\Bigl[\exp\Bigl\{\frac 1{\alpha_n^2}\sum_{k=0}^{n-1} f _n(S_k)\Bigr\}\1\{\supp(\ell_n)\subset B\}\Bigr]\\
&=\sum_{z_1,\dots,z_{n\alpha_n^{-2}}\in B}\prod_{i=1}^{n\alpha_n^{-2}}\E_{z_{i-1}}\Bigl[\exp\Bigl\{\frac 1{\alpha_n^2}\sum_{k=0}^{\alpha_n^2-1} f _n(S_k)\Bigr\}\1\{\supp(\ell_{\alpha_n^2})\subset B\}\1\{S_{\alpha_n^2}=z_i\}\Bigr]\\
&=\int_{Q_R^{n\alpha_n^{-2}}}\d x_1\dots\d x_{n\alpha_n^{-2}}\,\prod_{i=1}^{n\alpha_n^{-2}}\Bigl[\alpha_n^d
\E_{\lfloor x_{i-1}\alpha_n\rfloor}\Bigl[\exp\Bigl\{\frac 1{\alpha_n^2}\sum_{k=0}^{\alpha_n^2-1} f _n(S_k)\Bigr\}\\
&\qquad\qquad\qquad\qquad\times\1\{\supp(\ell_{\alpha_n^2})\subset B\}\1\{S_{\alpha_n^2}=\lfloor x_{i}\alpha_n\rfloor\}\Bigr)\Bigr].
\end{aligned}
\end{equation}
Let $(B_t)_{t\geq 0}$ be the Brownian motion on $\R^d$ with covariance matrix $\Gamma$, and let ${\tt E}_x$ denote the corresponding expectation, when $B_0=x\in\R^d$. Then $(\alpha_n^{-1}S_{\lfloor t\alpha_n^2\rfloor})_{t\geq 0}$ converges weakly towards $(B_t)_{t\geq 0}$ in distribution, and from a local central limit theorem (see \cite[P7.9, P7.10]{S76}) it follows that, uniformly in $x,y\in Q_R$,
\begin{equation}
\begin{aligned}
\lim_{n\to\infty}\alpha_n^d&\E_{\lfloor x\alpha_n\rfloor}\Bigl[\exp\Bigl\{\frac 1{\alpha_n^2}\sum_{k=0}^{\alpha_n^2-1} f _n(S_k)\Bigr\}\1\{\supp(\ell_{\alpha_n^2})\subset B\}\1\{S_{\alpha_n^2}=\lfloor y\alpha_n\rfloor\}\Bigr]\\
&={\tt E}_x\Bigl(\exp\Bigl\{\int_0^1 f (B_s)\,\d s\Bigr\}\1\{B_{[0,1]}\subset Q_R\};B_1\in \d y\Bigr)\Big/ \d y.
\end{aligned}
\end{equation}
Substituting this on the right hand side of \eqref{expansionexp} and again using the Markov property, we obtain that, as $n\to\infty$,
\begin{equation}\label{expansion2}
\E\Bigl[\exp\Bigl\{\frac 
n{\alpha_n^2} \langle f ,L_n\rangle\Bigr\}\1\{\supp(L_n)\subset Q_R\}\Bigr]=e^{o(n\alpha_n^{-2})}{\tt E}_0\Bigl(\exp\Bigl\{\int_0^{n\alpha_n^{-2}} f (B_s)\,\d s\Bigr\}\1\{B_{[0,n\alpha_n^{-2}]}\subset Q_R\}\Bigr).
\end{equation}
It is well-known that the expectation on the right hand side of \eqref{expansion2} is equal to $\exp\{\frac n{\alpha_n^2}[\lambda_R( f )+o(1)]\}$ as $n\to\infty$, and this ends the proof of Lemma~\ref{lem-cumul}.
\end{Proof}\qed

\begin{Proof}{Proof of Lemma~\ref{LDP}(i)} We shall 
apply a version of the abstract G\"artner-Ellis theorem (see 
\cite[Sect.~4.5]{DZ93}). (There is no problem in applying that result for {\em 
sub\/}probability measures instead of probability measure.) More precisely, 
we shall apply \cite[Cor.~4.5.27]{DZ93}, which implies the statement of 
Lemma~\ref{LDP}(i) under the following two assumptions: (1) the 
distributions of $L_n$ under $\P(\cdot\, ,\supp(L_n)\subset Q_R)$ form an exponentially tight 
family, and (2) the limit in \eqref{cumul} exists and is a finite, 
G\^ateau-differentiable and lower semicontinuous function of $ f $. These 
two points are satisfied in our case. Indeed, (1) is trivially satisfied 
since we consider subprobability measures on a compact set $Q_R$, and 
(2) follows from Lemma~\ref{lem-cumul}, together with \cite{Ga77}, where 
the G\^ateau-differentiability and lower semicontinuity of the map 
$ f \mapsto \lambda_R( f )$ is shown. An application of
\cite[Cor.~4.5.27]{DZ93} therefore yields the validity of a large deviation 
principle as stated in Lemma~\ref{LDP}(i).

It remains to identify the rate function obtained in 
\cite[Cor.~4.5.27]{DZ93} with the rate function of Lemma~\ref{LDP}(i), $\Ical_R$. 
The rate function appearing in  \cite[Cor.~4.5.27]{DZ93}, $\widetilde 
\Ical_R$, is the Legendre transform of $\lambda_R(\cdot)$:
\begin{equation}
\widetilde \Ical_R( \psi^2 )=\sup_{ f \in 
\Ccal(Q_R)}\bigl[\langle \psi^2 , f \rangle-\lambda_R( \psi^2 )\bigr],\qquad  \psi^2 \in\Fcal_R.
\end{equation}
It is obvious from \eqref{Lapeigenv} that $\lambda_R(\cdot)$ is itself 
the Legendre transform of $\Ical_R$, since $\Ical_R$ is equal to 
$\infty$ outside $\Fcal_R$. Because of the convexity inequality for gradients 
(see \cite[Theorem~7.8]{LL97}), $\Ical_R$ is a convex function on $\Fcal_R$. 
According to the Duality Lemma \cite[Lemma~4.5.8]{DZ93}, the Legendre 
transform of $\lambda_R(\cdot)$ is equal to $\Ical_R$, i.e., we have that 
$\widetilde \Ical_R= \Ical_R$. This finishes the proof of 
Lemma~\ref{LDP}(i).
\end{Proof}\qed

\begin{Proof}{Proof of Lemma~\ref{LDP}(ii)}
This is a modification of the proof of part (i) above; we point out the 
differences only.  Recall that we identify the 
box $B_R=\{\lfloor -R\rfloor+1,\dots,\lfloor R\rfloor-1\}^d$ with the 
torus $\{\lfloor -R\rfloor+1,\dots,\lfloor R\rfloor\}^d$ where $\lfloor 
R\rfloor$ is identified with $\lfloor -R\rfloor+1$. Analogously, we
conceive $Q_R=[-R,R]^d$ as the $d$-dimensional torus with the opposite
sides identified.

For a continuous bounded function $ f \colon Q_R\to\R$, 
introduce the principal eigenvalue of the operator 
$\frac{1}{2}\nabla\cdot\Gamma\nabla  + f $ on $L^2(Q_R)$ with periodic boundary 
condition:
\begin{equation}
\lambda^{\smallsup{R}}( f )=\max\Bigl\{\int_{Q_R} f (x)
\psi^2(x)\,\d x-\frac{1}{2}\int_{Q_R}\big|\Gamma^{\frac 12}\nabla_R \psi(x)\big|^2\,\d x\colon  \psi\in\Ccal_1(Q_R),\int_{Q_R}\psi^2(x)\,\d x=1\Bigr\},
\end{equation}
where we recall that $\nabla_R$ is the gradient of the torus $Q_R$.

The main step in the proof of Lemma~\ref{LDP}(ii) is to show 
that, for any continuous bounded function $ f \colon  
Q_R\to\R$,
\begin{equation}\label{eigenasyR}
\lambda^{\smallsup{R}}( f )=\lim_{n\to\infty}\frac 
{\alpha_n^2}n\log \E\Bigl[\exp\Bigl\{\frac 
n{\alpha_n^2}\langle f ,L_n^{\smallsup{R\alpha_n}}\rangle\Bigr\}\Bigr].
\end{equation}
This is done in the same way as in the proof of Lemma~\ref{lem-cumul}, noting that 
the process $(\alpha_n^{-1}S^{\smallsup{R\alpha_n}}_{t\alpha_n^2})_{t\geq 0}$ converges weakly in distribution towards $(B_t^{\smallsup{R}})_{t\ge 0}$, the Brownian motion with covariance matrix $\Gamma$, wrapped around the torus $Q_R$. Also using a local central limit theorem, we obtain, as $n\to\infty$,
\begin{equation}
\E\Bigl[\exp\Bigl\{\frac 
n{\alpha_n^2}\langle f ,L_n^{\smallsup{R\alpha_n}}\rangle\Bigr\}\Bigr)=e^{o(n\alpha_n^{-2})}
{\tt E_0}\Bigl(\exp\Bigl\{\int_0^{n\alpha_n^{-2}} f (B_s^{\ssup{R}})\,\d s\Bigr\}\Bigr].
\end{equation}
It is well-known that the expectation on the right side is equal to $\exp\{\frac n{\alpha_n^2}[\lambda^{\smallsup{R}}( f )+o(1)]\}$ as $n\to\infty$, and this shows that also \eqref{eigenasyR} holds. The remainder of the proof of Lemma~\ref{LDP}(ii) is the same as the proof of Lemma~\ref{LDP}(i).
\end{Proof}\qed

\noindent {\bf Acknowledgment. } This work was partially supported by DFG grant Ko 2205/1-1. W.~K. thanks the Deutsche Forschungsgemeinschaft for awarding a Heisenberg grant (realized in 2003/04). W.~K.\ and N.~G.\ thank the Laboratoire de Probabilit{\'e}s for its hospitality. All three authors thank Francis Comets for helpful discussions. We also thank the referee for 
carefully reading the first version of the paper.

\end{document}